\documentclass[letterpaper,11pt]{article}
\usepackage{amsmath}
\usepackage{amssymb}
\usepackage{geometry}
\usepackage[latin1]{inputenc}
\usepackage[T1]{fontenc}
\usepackage{times}
\usepackage[colorlinks=true,linkcolor=red,citecolor=blue]{hyperref}
\usepackage{natbib}
\usepackage{bbm}
\usepackage{enumitem}

\date{}
\newtheorem{theorem}{Theorem}[section]

\newtheorem{corollary}[theorem]{Corollary}

\newtheorem{remark}[theorem]{Remark}
\newtheorem{theo}{Theorem}

\numberwithin{equation}{section}
\newcommand{\lqn}[1]{\noalign{\noindent $\displaystyle#1$}}

\begin{document}

\title{Some asymptotic results for the integrated empirical process
with applications to statistical tests}
\author{Sergio Alvarez-Andrade\footnote{e-mail: sergio.alvarez@utc.fr}~,
~ Salim Bouzebda\footnote{e-mail: salim.bouzebda@utc.fr}
~ and Aim\'e Lachal\footnote{e-mail: aime.lachal@insa-lyon.fr}
\\~\\
$^{*\dag}$ Sorbonne Universit\'es, Universit\'e de Technologie de Compi\`egne
\\
\textit{Laboratoire de Math\'ematiques Appliqu\'ees de Compi\`egne}
\\
$^\ddag$ Universit\'e de Lyon, Institut National des Sciences Appliquées de Lyon
\\
\textit{Institut Camille Jordan}}
\date{}
\maketitle

\begin{abstract}
\noindent
The main purpose of this paper is to investigate the strong approximation of
the integrated empirical process. More precisely, we obtain the exact rate of the
approximations by a sequence of weighted Brownian bridges and a weighted Kiefer
process. Our arguments are based in part on the \cite{KMT1975}'s results.
Applications include the two-sample testing procedures together with the
change-point problems. We also consider the strong approximation of the integrated
empirical process when the parameters are estimated. Finally, we study the
behavior of the self-intersection local time of the partial sum process
representation of the integrated empirical process.
\\
\noindent \textbf{Key words:} Integrated empirical process; Brownian bridge;
Kiefer process; Rates of convergence; Local time; Two-sample problem;
Hypothesis testing; Goodness-of-fit; Change-point.
\\
\noindent \textbf{AMS Classifications:} primary: 62G30; 62G20; 60F17;
secondary: 62F03; 62F12; 60F15.
\end{abstract}

\section{Introduction}

Let $\{X_i: i\in\mathbb{N}^*\}$ be a sequence of independent, identically
distributed [i.i.d.] random variables [r.v.'s] with continuous distribution
function [d.f.] $F(t)=\mathbb{P}\{X_1\leq t\}$ for $t\in\mathbb{R}$. For each
$n\in\mathbb{N}^*$, let $\mathbb{F}_n$ be the empirical d.f. based upon the
sample $X_1,\ldots,X_n$ defined as
\begin{equation}
\mathbb{F}_n(t) :=\frac{1}{n}\#\{i\in\{1,\dots,n\}:X_i\leq t\}
=\frac{1}{n}\sum_{i=1}^n \mathbbm{1}_{\{X_i\leq t\}}
\quad\text{for}\quad t\in\mathbb{R},
\end{equation}
where $\#$ denotes cardinality. For each $n\in\mathbb{N}^*$,
we introduce the \emph{empirical process} $\alpha_n$ defined by
\begin{equation}\label{beta1}
\alpha_n(t):=\sqrt{n}\,(\mathbb{F}_n(t)-F(t))
\quad\text{for}\quad t\in\mathbb{R}.
\end{equation}
We denote by $Q(y) :=\inf\{x\in\mathbb{R}: F(x)\geq y\}$ for $y\in(0,1)$
the usual quantile function (generalized inverse) pertaining to $F$, and extend its
definition to $[0,1]$ by setting
$$
Q(0):=\lim_{y\downarrow 0}Q(y)\quad\text{and}\quad Q(1):=\lim_{y\uparrow 1}Q(y).
$$
Consider additionally a sequence of i.i.d. r.v.'s $\{U_i: i\in\mathbb{N}^*\}$
uniformly distributed on $[0,1]$ such that $X_i = Q(U_i)$ (cf.
\cite{Shorack1986}, p.~3 and the references therein). For each
$n\in\mathbb{N}^*$, let $\mathbb{U}_n$ be the empirical d.f. based upon the
sample $U_1,\ldots,U_n$:
$$
\mathbb{U}_n(u) := \frac{1}{n}\#\{i\in\{1,\dots,n\}:U_i \leq u\}
=\frac{1}{n}\sum_{i=1}^n \mathbbm{1}_{\{U_i\leq u\}}
\quad\text{for}\quad u\in[0,1].
$$
We introduce the corresponding \emph{uniform empirical process} $\beta_n$ defined,
for each $n\in\mathbb{N}^*$, by
\begin{equation}\label{uniformprocesse}
\beta_n(u) := \sqrt{n}\left(\mathbb{U}_n(u)-u\right)
\quad\text{for}\quad u\in[0,1].
\end{equation}
We have the usual relation between the empirical process and uniform empirical
process:
$$
\alpha_n(t)=\beta_n(F(t))\quad\text{for}\quad t\in \mathbb{R},\,n\in\mathbb{N}^*.
$$

In this paper, we consider the \textit{integrated empirical d.f.} based upon
the sample $X_1,\ldots,X_n$
\begin{equation}\label{beta2}
\overline{\mathbb{F}}_n(t):=\int_{-\infty}^t\mathbb{F}_n(s)\,d\mathbb{F}_n(s)
\quad\text{for}\quad t\in \mathbb{R},\,n\in\mathbb{N}^*
\end{equation}
together with the corresponding \textit{integrated empirical process}
\begin{equation}\label{beta12}
\overline{\alpha}_n(t):=\sqrt{n}\left(\overline{\mathbb{F}}_n(t)-\overline{F}(t)\right)
\quad\text{for}\quad t\in \mathbb{R},\,n\in\mathbb{N}^*
\end{equation}
where
$$
\overline{F}(t):=\int_{-\infty}^t F(s)\,dF(s)=\frac{1}{2}\,F(t)^2.
$$
From Theorem II.1 of \cite{HenzeNikitin2002} (see also the self-contained proof
given in Section~\ref{appe}), we learn that, almost surely,
\begin{equation}\label{beta21}
\overline{\mathbb{F}}_n(t)=\frac{1}{2}\left(\mathbb{F}_n(t)^2
+\frac{1}{n}\,\mathbb{F}_n(t) \right)
\quad\text{for}\quad t\in \mathbb{R},\,n\in\mathbb{N}^*.
\end{equation}

\cite{HenzeNikitin2000,HenzeNikitin2002} introduced and deeply investigated
the goodness-of-fit testing procedures based on the integrated empirical process.
Indeed, the asymptotic properties of their procedures, Kolmogorov-Smirnov,
Cram\'er-von Mises and Watson-type statistics, can be derived from the limiting
behavior of the integrated empirical process. \cite{HenzeNikitin2003} considered
a two-sample testing procedure and focused on the approximate local Bahadur
efficiencies of their statistical tests. It is noteworthy to point out that
tests based on some integrated empirical processes turn out to be more efficient
for certain distributions. In \cite{Aime2001}, another version of the
integrated empirical process (actually a broader class of $p$-fold integrated
empirical process, $p\in\mathbb{N}^*$) was introduced. For the extension
to the multivariate framework, we may refer to \cite{Jing2006} and \cite{Jing2007}
where some projected integrated empirical processes for testing the equality of
two multivariate distributions are considered. Inspired by the work of
\cite{HenzeNikitin2003}, \cite{BouzebdaElfaouzy2012} developed multivariate
two-sample testing procedures based on the integrated empirical copula process
that are extended to the $K$-sample problem in \cite{Bouzebdaetal2011}.
Emphasis is placed on the explanation of the strong approximation methodology.
The asymptotic behavior of weighted multivariate Cram\'er-von Mises-type
statistics under contiguous alternatives was characterized by
\cite{BouzebdaZari2013}.

The main purpose of this paper is to investigate the strong approximation of
the integrated empirical process. Next we study the asymptotic properties
of statistical tests based on this process. We also study the behavior of
the self-intersection local time of the partial sum process representation of the
integrated uniform empirical process. To our best knowledge, the problem that
we consider was open up to presently, and it gives the main motivation to our study.

Let us point out that the strong approximations are quite useful and have
received considerable attention in probability theory. Indeed, many
well-known and important probability theorems can be considered as consequences of
results about strong approximation of sequences of sums by corresponding
Gaussian sequences.

We will first obtain an upper bound in probability for the distance
between the integrated empirical process and a sequence of appropriate
Brownian bridges (see Theorem~\ref{lem2}). This is the key point of our study.
From this, we will deduce a strong approximation of the integrated
empirical process by this sequence of Brownian bridges (see Corollary~\ref{corollaryapprox}).
As an application, we will derive the rates of convergence for the
distribution of smooth functionals of the integrated empirical process
(see Corollary~\ref{corol22}). Moreover, we will deduce
strong approximations for the Kolmogorov-Smirnov and Cram\'er-von Mises-type
statistics associated with the integrated empirical process
(see Corollary~\ref{ecooooo}).

Second, we will obtain a strong approximation of the integrated
empirical process by a Kiefer processes (see Theorem~\ref{kieferapproximation}).
This latter is of particular interest; indeed, for instance, any kind of law
of the iterated logarithm which holds for the partial sums of Gaussian
processes may then be transferred to the integrated empirical process
(see Corollary \ref{th2}).
We may refer to \cite{Dasgupta2008} (Chapter 12), \cite{Csorgoho1993} (Chapter 3),
\cite{Csorgo1981REVESZ} (Chapters 4-5) and \cite{Shorack1986} (Chapter 12) for
expositions, details and references about this problem.

We refer to \cite{Hall1983S}, \cite{csorgo2007} and \cite{Mason2012} for a
survey of some applications of the strong approximation and many references.
There is a huge literature on the strong approximations and their applications.
It is not the purpose of this paper to survey this extensive literature.

The layout of the article is as follows. In Section~\ref{sectionone}, we first
present some strong approximation results for the integrated empirical
process; our main tools are the results of \cite{KMT1975}. Sections~\ref{sectiontwo}
and \ref{sectionchnae} are devoted to statistical applications,
namely the two-sample and change-point problems respectively. In Section~\ref{Section5},
we deal with the strong approximation of the integrated
empirical process when parameters are estimated. Section~\ref{section343} is
concerned with the behavior of the self-intersection local time of the partial
sum process representation of the integrated empirical process.
To prevent from interrupting the flow of the presentation, all mathematical
developments are postponed to Section~\ref{appe}.

\section{Strong approximation}\label{sectionone}

First, we introduce some definitions and notations. Let $\mathbb{W}=
\{\mathbb{W}(s): s\geq 0\}$ and $\mathbb{B}=\{\mathbb{B}(u) : u\in[0,1]\}$ be
the standard Wiener process and Brownian bridge, that is, the centered Gaussian
processes with continuous sample paths, and covariance functions
$$
\mathbb{E}(\mathbb{W}(s)\mathbb{W}(t))=s\wedge t\quad\text{for}\quad s,t\geq 0
$$
and
$$
\mathbb{E}(\mathbb{B}(u)\mathbb{B}(v))=u\wedge v- uv\quad\text{for}\quad u,v\in[0,1].
$$
A Kiefer process $\mathbb{K}=\{\mathbb{K}(s,u) : s \geq 0, u\in[0,1]\}$ is a two-parameters
centered Gaussian process, with continuous sample paths, and covariance function
$$
\mathbb{E}(\mathbb{K}(s,u)\mathbb{K}(t,v))=(s\wedge t)\,(u\wedge v-uv)
\quad\text{for}\quad s,t\geq 0\quad\text{and}\quad u,v\in[0,1].
$$
It satisfies the following distributional identities:
$$
\left\{\mathbb{K}(s,u): u\in[0,1]\right\}\stackrel{\mathcal{L}}{=}
\left\{\sqrt{s}\,\mathbb{B}(u): u\in[0,1]\right\}\quad\text{for}\quad s \geq 0
$$
and
$$
\left\{\mathbb{K}(s,u): s\geq 0\right\}\stackrel{\mathcal{L}}{=}
\left\{\sqrt{u(1-u)}\,\mathbb{W}(s): s\geq 0\right\}\quad\text{for}\quad u\in[0,1],
$$
where $\stackrel{\mathcal{L}}{=}$ stands for the equality in distribution.
The interested reader may refer to \cite{Csorgo1981REVESZ} for details on
the Gaussian processes mentioned above. It is well-known that the empirical
uniform process $\{\beta_n:n\in\mathbb{N}^*\}$ converges to $\mathbb{B}$ in
$D[0,1]$ (the space of all right-continuous real-valued functions defined on
$[0,1]$ which have left-hand limits, equipped with the Skorohod topology; see,
for details, \cite{Billingsley1968}). The rate of convergence of this process
to $\mathbb{B}$ is an important task in statistics as well as in probability
that has been investigated by several authors. We can and will assume without
loss of generality that all r.v.'s and processes introduced so far and later on
in this paper can be defined on the same probability space (cf. Appendix 2 in
\cite{Csorgoho1993}). Koml\'os, Major, and Tusn\'ady [KMT] (\cite{KMT1975},
Theorem 3; refer also to \cite{Major1976}) stated the following Brownian bridge
approximation for $\{\beta_n:n\in\mathbb{N}^*\}$, along with a description
of its proof with few details.

\begin{theo}\label{theoA}
On a suitable probability space, we can define the uniform empirical process
$\{\beta_n:n\in\mathbb{N}^*\}$, in combination with a sequence of Brownian
bridges $\left\{\mathbb{B}_n:n\in\mathbb{N}^*\right\}$, such that,
for all $x>0$ and all $n\in\mathbb{N}^*$,
\begin{equation}\label{major}
\mathbb{P}\bigg\{\sup_{u\in[0,1]} |\beta_n(u)-\mathbb{B}_n(u)|
\geq \frac{1}{\sqrt{n}}\,(c_1\log n+x)\bigg\}\leq c_2 \,e^{-c_3x},
\end{equation}
where $c_1$, $c_2$ and $c_3$ are positive universal constants.
\end{theo}
In his manuscript, \cite{Major2000} details the original proof of (\ref{major}).
\cite{Chatterjee2012} provided a new alternative approach for proving the
famous KMT theorem.

\begin{remark}
In the sequel, the precise meaning of ``suitable probability space'' is that
an independent sequence of Wiener processes, which is independent of the
originally given sequence of i.i.d. r.v.'s, can be constructed on the assumed
probability space. This is a technical requirement which allows the construction
of the Gaussian processes displayed in our theorems, and which is not restrictive
since one can expand the probability space to make it rich enough (see, e.g.,
Appendix 2 in \cite{Csorgoho1993}, \cite{deAcosta1982}, \cite{Csorgo1981REVESZ}
and Lemma A1 in \cite{Berkes1979}). Throughout this paper, it will be assumed
that the underlying probability spaces are suitable in this sense.
\end{remark}

In the following theorem, we state the key point to access the strong
approximation of the integrated empirical process
$\{\overline{\alpha}_n,n\in\mathbb{N}^*\}$.
\begin{theorem}\label{lem2}
On a suitable probability space, we may define the
integrated empirical process $\{\overline{\alpha}_n:n\in\mathbb{N}^*\}$,
in combination with a sequence of Brownian bridges
$\left\{\mathbb{B}_n:n\in\mathbb{N}^*\right\}$, such that,
for large enough $x$ and all $n\in\mathbb{N}^*$,
\begin{equation}\label{estimation}
\mathbb{P}\!\left\{ \sup_{t\in \mathbb{R}} \big|\overline{\alpha}_n(t)
-F(t)\,\mathbb{B}_n(F(t))\big| \geq
\frac{1}{\sqrt{n}}\,(A\log n+x)\right\}
\leq B\,e^{-Cx}
\end{equation}
where $A$, $B$ and $C$ are positive universal constants.
\end{theorem}

An important consequence of Theorem~\ref{lem2} is an upper bound for the
convergence of distributions of smooth functionals of $\overline{\alpha}_n$.
Indeed, applying (\ref{estimation}) with $x=c\,\log n$
for a suitable constant $c$ yields the result below.
\begin{corollary}\label{corol22}
If $\Phi(\cdot)$ is a Lipschitz functional defined
on  $D[0,+\infty)$ such that the r.v. $\Phi(F(\cdot)\,\mathbb{B}(F(\cdot)))$
admits a bounded density function, then, as $n\to\infty$,
\begin{equation}\label{estimationfunctional}
\sup_{x\in\mathbb{R}}\left|\mathbb{P}\!\left\{\Phi\big(\overline{\alpha}_n(\cdot)\big)
\leq x\right\}-\mathbb{P}\!\left\{\Phi\big(F(\cdot)\,
\mathbb{B}(F(\cdot))\big)\leq x\right\}\right|
=O\!\left(\frac{\log n}{\sqrt{n}}\right)\!.
\end{equation}
\end{corollary}
For more comments on this kind of results, we may refer to \cite{Kokoszka2000},
Corollary 1.1 and p.~2459.

By applying (\ref{estimation}) to $x=c'\log n$ for a suitable constant $c'$
and appealing to Borel-Cantelli lemma, one can obtain the following almost sure
approximation of the process $\{\overline{\alpha}_n:n\in\mathbb{N}^* \}$ based on a
sequence of Brownian bridges.
\begin{corollary}\label{corollaryapprox}
The following bound holds, with probability~$1$, as $n\to\infty$:
\begin{equation}\label{approx}
\sup_{t\in \mathbb{R}} \big|\overline{\alpha}_n(t)- F(t)\,\mathbb{B}_n(F(t))\big|
=O\!\left(\frac{\log n}{\sqrt{n}}\right)\!.
\end{equation}
\end{corollary}

The next result yields an almost sure approximation for
$\{\overline{\alpha}_n:n\in\mathbb{N}^*\}$ based on a Kiefer process.
\begin{theorem}\label{kieferapproximation}
On a suitable probability space, we may define the integrated empirical process
$\{\overline{\alpha}_n:n\in\mathbb{N}^*\}$, in combination with
a Kiefer process $\{\mathbb{K}(s,u): s\geq 0, u\in[0,1]\}$, such that,
with probability~$1$, as $n\to\infty$,
$$
\max_{1\leq k\leq n}\sup_{t\in \mathbb{R}}\left|\sqrt{k}\,\overline{\alpha}_k(t)
-F(t)\,\mathbb{K}(k,F(t))\right|=O\!\left((\log n)^2\right)\!.
$$
\end{theorem}
From Theorem~\ref{kieferapproximation} and by invoking the law of
the iterated logarithm for Gaussian sequences, we have almost surely
\begin{align}
\limsup_{n\to\infty}\frac{\sup_{t\in\mathbb{R}}
\big|\overline{\alpha}_n(t)\big|}{\sqrt{\log \log n}}
&=\limsup_{n\to\infty}\frac{\sup_{t\in \mathbb{R}}
\big|F(t)\,\mathbb{K}(n,F(t))\big|}{\sqrt{n\log \log n}}
\nonumber\\
&=\sup_{t\in \mathbb{R}} \sqrt{2\,{\rm Var}\big(F(t)\,\mathbb{K}(1,F(t))\big)}
=\sup_{u\in[0,1]} \sqrt{2\,{\rm Var}\big(u\,\mathbb{K}(1,u)\big)}.
\label{lil}
\end{align}
Observing that ${\rm Var}(u\,\mathbb{K}(1,u))=u^3(1-u)$ and
$\sup_{u\in[0,1]}u^3(1-u)=27/256$, (\ref{lil}) readily implies the following corollary
(``a.s.'' stands for ``almost surely'').
\begin{corollary}\label{th2}
We have the following law of iterated logarithm for the integrated
empirical process:
\begin{equation}\label{lli}
\limsup_{n\to\infty}\frac{\sup_{t\in\mathbb{R}}\big|\overline{\alpha}_n(t)\big|}{\sqrt{\log \log n}}
=\frac{3\sqrt{3}}{8\sqrt{2}}\quad\text{a.s.}
\end{equation}
\end{corollary}

As a direct application of (\ref{approx}) and (\ref{lli}) to the problem of
goodness-of-fit, for testing the null hypothesis
$$
\mathcal{H}_0: F=F_0,
$$
we can use the following statistics: the \emph{integrated
Kolmogorov-Smirnov statistic} as well as the \emph{integrated
Cram\'er-von Mises statistic}
$$
\overline{\mathbf{S}}_n:=\sup_{t\in \mathbb{R}}\left|\sqrt{n}\left(\overline{\mathbb{F}}_n(t)
-\overline{F}_{0}(t)\right)\right|\quad\text{and}\quad
\overline{\mathbf{T}}_n:=n\int_{\mathbb{R}} \left(\overline{\mathbb{F}}_n(t)
-\overline{F}_{0}(t)\right)^2dF_{0}(t).
$$
\begin{corollary}\label{ecooooo}
We have, under $\mathcal{H}_0$, with probability~$1$, as $n\to\infty$,
\begin{align}
\left|\overline{\mathbf{S}}_n-\sup_{t\in\mathbb{R}}
\big|F(t)\,\mathbb{B}_n(F(t))\big|\right|
&
=O\!\left(\frac{\log n}{\sqrt{n}}\right)\!,
\label{kol}\\[1ex]
\left|\overline{\mathbf{T}}_n-\int_{\mathbb{R}}
\big[F(t)\,\mathbb{B}_n(F(t))\big]^2\,dF_{0}(t)\right|
&
=O\!\left(\!\sqrt{\frac{\log\log n}{n}}\,\log n\right)\!.
\label{vonmises}
\end{align}
\end{corollary}

\section{The two-sample problem}\label{sectiontwo}

For each $m,n\in\mathbb{N}^*$, let $X_1,\ldots, X_m$ and $Y_1,\ldots,Y_n$
be independent random samples from continuous d.f.'s $F$ and $G$, respectively,
and let $\overline{\mathbb{F}}_m$ and $\overline{\mathbb{G}}_n$ denote their respective
integrated empirical d.f.'s. Tests for the null hypothesis
$$
\mathcal{H}_0': F=G,
$$
can be based on the \emph{integrated two-sample empirical process}
defined, for each $m,n\in\mathbb{N}^*$, by
$$
\overline{\boldsymbol{\xi}}_{m,n}(t):=\sqrt{\frac{mn}{m+n}}
\left(\overline{\mathbb{F}}_n(t)-\overline{\mathbb{G}}_m(t)\right)
\quad\text{for}\quad t\in \mathbb{R}.
$$
We can use the following statistics for testing $\mathcal{H}_0'$: the
\emph{integrated two-sample Kolmogorov-Smirnov statistic}
as well as the \emph{integrated two-sample Cram\'er-von Mises statistic}
$$
\overline{\mathbf{S}}_{m,n}:=\sup_{t\in \mathbb{R}}
\left|\overline{\boldsymbol{\xi}}_{m,n}(t)\right|\quad\text{and}\quad
\overline{\mathbf{T}}_{m,n}:=\int_{\mathbb{R}}\overline{\boldsymbol{\xi}}_{m,n}(t)^2\,dF(t).
$$
Set, for any $m,n\in\mathbb{N}^*$,
$$
\varphi(m,n):=\max\left(\frac{\log m}{\sqrt{m}},\frac{\log n}{\sqrt{n}}\right)
\quad\text{and}\quad\phi(m,n):=\max\!\left(\sqrt{\frac{\log\log m}{m}}\,\log m,
\sqrt{\frac{\log\log n}{n}}\,\log n\right)\!.
$$

The following results are consequences of Corollary~\ref{corollaryapprox}.
\begin{corollary}\label{ximn-approx}
On a suitable probability space, it is possible to define
$\big\{\overline{\boldsymbol{\xi}}_{m,n}:m,n\in\mathbb{N}^*\big\}$,
jointly with a sequence of Gaussian processes $\{\mathbb{B}_{m,n}^*:m,n\in\mathbb{N}^*\}$,
such that, under $\mathcal{H}_0'$, with probability $1$, as $\min(m,n) \to\infty$,
$$
\sup_{t\in\mathbb{R}}\left|\overline{\boldsymbol{\xi}}_{m,n}(t)-\mathbb{B}_{m,n}^*(t)\right|=O(\varphi(m,n)),
$$
where
$$
\mathbb{B}_{m,n}^*(t)=F(t)\left(\sqrt{\frac{n}{n+m}}\,\mathbb{B}_m^1(F(t))-\sqrt{\frac{m}{n+m}}\,\mathbb{B}_n^2(F(t))\right)\!,
$$
the processes $\{\mathbb{B}_m^1:m\in\mathbb{N}^*\}$ and $\{\mathbb{B}_n^2:n\in\mathbb{N}^*\}$
consisting of two sequences of Brownian bridges constructed as in Theorem~\ref{lem2}.
\end{corollary}
\begin{corollary}\label{ximn-stat}
We have, under $\mathcal{H}_0'$, with probability $1$, as $\min(m,n) \to\infty$,
$$
\left|\overline{\mathbf{S}}_{m,n}-\sup_{t\in \mathbb{R}}
|\mathbb{B}_{m,n}^*(t)|\right|=O(\varphi(m,n))\quad\text{and}\quad
\left|\overline{\mathbf{T}}_{m,n}-\int_{\mathbb{R}}
\mathbb{B}_{m,n}^*(t)^2\,dF(t)\right|=O(\phi(m,n)).
$$
\end{corollary}

As in \cite{BouzebdaElfaouzy2012}, consider the following \emph{modified integrated two-sample
empirical process}, for a fixed positive integer $q$,
$$
\overline{\boldsymbol{\xi}}_{m,n}^{(q)}(t):=\sqrt{\frac{mn}{m+n}}
\left(\overline{\mathbb{F}}_m(t)^q-\overline{\mathbb{G}}_n(t)^q\right)
\quad\text{for}\quad t\in \mathbb{R}.
$$
Reasonable statistics for testing $\mathcal{H}_0'$ would be the
\emph{modified integrated Kolmogorov-Smirnov statistic}
and the \emph{modified integrated Cram\'er-von Mises statistic}
$$
\overline{\mathbf{S}}_{m,n}^{(q)}:=
\sup_{t\in \mathbb{R}}\left|\overline{\boldsymbol{\xi}}_{m,n}^{(q)}(t)\right|
\quad\text{and}\quad
\overline{\mathbf{T}}_{m,n}^{(q)}:=\int_{\mathbb{R}}
\overline{\boldsymbol{\xi}}_{m,n}^{(q)}(t)^2\,dF(t).
$$

We extend Corollary~\ref{ximn-approx} and \ref{ximn-stat} as follows.
\begin{corollary}\label{ximn}
On a suitable probability space, it is possible to define
$\big\{\overline{\boldsymbol{\xi}}_{m,n}^{(q)}:m,n\in\mathbb{N}^*\big\}$,
jointly with a sequence of Gaussian processes $\big\{\mathbb{B}_{m,n}^{*(q)}:m,n\in\mathbb{N}^*\big\}$,
such that, under $\mathcal{H}_0'$, with probability $1$, as $\min(m,n)\to\infty$,
$$
\sup_{t\in\mathbb{R}}\left|\overline{\boldsymbol{\xi}}_{m,n}^{(q)}(t)-\mathbb{B}_{m,n}^{*(q)}(t)\right|=O(\varphi(m,n)),
$$
where
$$
\mathbb{B}_{m,n}^{*(q)}(t):=\frac{q}{2^{q-1}}\,F(t)^{2q-1}\left(\sqrt{\frac{n}{m+n}}\,
\mathbb{B}_m^1(F(t))-\sqrt{\frac{m}{m+n}}\,\mathbb{B}_n^2(F(t))\right)\!.
$$
\end{corollary}
\begin{corollary}
Under $\mathcal{H}_0'$, with probability $1$, as $\min(m,n) \to\infty$, we have
$$
\left|\overline{\mathbf{S}}_{m,n}^{(q)}-\sup_{t\in\mathbb{R}}
\big|\mathbb{B}_{m,n}^{*(q)}(t)\big|\right|=O(\varphi(m,n))\quad\text{and}\quad
\left|\overline{\mathbf{T}}_{m,n}^{(q)}-\int_{\mathbb{R}}
\mathbb{B}_{m,n}^{*(q)}(t)^2\,dF(t)\right|=O(\phi(m,n)).
$$
\end{corollary}
\begin{remark}
The family of statistics indexed by $q$ may be used to maximize the power of the
statistical test for a specific alternative hypothesis as argued in \cite{Ahmad2000}.
\end{remark}

Now, we fix a positive integer $K$ and we describe the more general $K$-sample problem.
For each $k\in\{1,\ldots,K\}$, we consider a setting made of independent observations
$\big\{X_i^k:i\in\{1,\ldots,n_k\}\big\}$ of a real-valued r.v. $X^k$.
The d.f. of $X^k$ is denoted by $F^k$ and is assumed to be continuous.
We would like to test, $F_0$ being a fixed continuous d.f., the null hypothesis
$$
\mathcal{H}_0^K: F^1= F^2=\dots=F^K=F_0.
$$
For any $K$-tuple of positive integers $\boldsymbol{n}=(n_1,\dots,n_K)$, set
$\boldsymbol{|n|}=\sum_{k=1}^K n_k$ and let
$$
(Z_1,\ldots,Z_{\boldsymbol{|n|}}):=\big(X_1^1,\ldots,X_{n_1}^1,X_1^2,\ldots,
X_{n_2}^2,\ldots,X_1^K,\ldots,X_{n_K}^K\big)
$$
be the pooled sample of total size $\boldsymbol{|n|}$,
$\overline{\mathbb{D}}_{K,\boldsymbol{n}}$ be the integrated empirical
d.f. based upon $Z_1,\ldots,Z_{\boldsymbol{|n|}}$, and, for each
$k\in\{1,\ldots,K\}$, $\overline{\mathbb{F}}_{n_k}^k$ be the integrated
empirical d.f. based upon $X_1^k,\ldots,X_{n_k}^k$. Of course,
we have the following identity:
\begin{equation}\label{DD}
\overline{\mathbb{D}}_{K,\boldsymbol{n}}=\frac{1}{\boldsymbol{|n|}}
\sum_{k=1}^K n_k\,\overline{\mathbb{F}}_{n_k}^k.
\end{equation}

Next, we define the \emph{integrated $K$-sample empirical process} in the
following way: for any $K$-tuple $\boldsymbol{n}=(n_1,\dots,n_K)\in(\mathbb{N}^*)^K$,
$$
\overline{\boldsymbol{\xi}}_{K,\boldsymbol{n}}(t):=\sum_{k=1}^Kn_k
\left(\overline{\mathbb{F}}_{n_k}(t)-\overline{\mathbb{D}}_{K,\boldsymbol{n}}(t)\right)^2
\quad\text{for}\quad t\in \mathbb{R}.
$$
Obvious candidates for testing Hypothesis $\mathcal{H}_0^K$ are the
\emph{integrated $K$-sample Kolmogorov-Smirnov statistics}
and the \emph{integrated $K$-sample Cram\'er-von Mises functional}
(the usual square being included in the definition of
$\overline{\boldsymbol{\xi}}_{K,\boldsymbol{n}}$)
$$
\overline{\mathbf{S}}_{K,\boldsymbol{n}}:=\sup_{t\in\mathbb{R}}
\overline{\boldsymbol{\xi}}_{K,\boldsymbol{n}}(t)\quad\text{and}\quad
\overline{\mathbf{T}}_{K,\boldsymbol{n}}:=\int_{\mathbb{R}}
\overline{\boldsymbol{\xi}}_{K,\boldsymbol{n}}(t)\,dF_0(t).
$$
Set
$$
\phi_K^*(\boldsymbol{n}):=\max_{1\leq k\leq K}
\left\{\sqrt{\frac{\log\log n_k}{n_k}}\,\log n_k \right\}\!.
$$

The following result is a consequence of Corollary~\ref{corollaryapprox}.
\begin{theorem}\label{theoremK}
On a suitable probability space, it is possible to define
$\big\{\overline{\boldsymbol{\xi}}_{K,\boldsymbol{n}}:\boldsymbol{n}\in(\mathbb{N}^*)^K\big\}$,
jointly with a sequence of processes $\big\{\mathbb{B}_{K,\boldsymbol{n}}^*:
\boldsymbol{n}\in(\mathbb{N}^*)^K\big\}$, such that, under $\mathcal{H}_0^K$, with
probability~$1$, as $\min_{1\leq k\leq K}n_k\to\infty$,
$$
\sup_{t\in\mathbb{R}}\left|\overline{\boldsymbol{\xi}}_{K,\boldsymbol{n}}(t)-\mathbb{B}_{K,\boldsymbol{n}}^*(t)\right|
=O\big(\phi_K^*(\boldsymbol{n})\big)\!,
$$
where, for each $\boldsymbol{n}=(n_1,\dots,n_K)\in(\mathbb{N}^*)^K$,
$\mathbb{B}_{K,\boldsymbol{n}}^*$ is the process defined by
$$
\mathbb{B}_{K,\boldsymbol{n}}^*(t):=F_0(t)^2\Bigg[\left(\sum_{k=1}^K
\mathbb{B}_{n_k}^k\!\big(F_0(t)\big)\right)^{\!\!2}-
\left(\sum_{k=1}^K\sqrt{\frac{n_k}{|\boldsymbol{n}|}}\,
\mathbb{B}_{n_k}^k\!\big(F_0(t)\big)\right)^{\!\!2}\Bigg]
\quad\text{for}\quad t\in\mathbb{R},
$$
the processes $\big\{\mathbb{B}_m^k:m\in\mathbb{N}^*\big\}$, $k\in\{1,\dots,K\}$,
consisting of $K$ sequences of Brownian bridges constructed as in Theorem~\ref{lem2}.
\end{theorem}

The next result, which is an immediate consequence of the previous theorem
(observe that $\overline{\mathbf{S}}_{K,\boldsymbol{n}}$ and
$\overline{\mathbf{T}}_{K,\boldsymbol{n}}$ are bounded linear functionals of the
process $\overline{\boldsymbol{\xi}}_{K,\boldsymbol{n}}$), gives the limit
null distributions of the statistics under consideration.
\begin{corollary}
We have, under $\mathcal{H}_0^K$, with probability~$1$, for
$\boldsymbol{n}=(n_1,\dots,n_K)$ such that $\min_{1\leq k \leq K} n_k \to\infty$,
$$
\left|\overline{\mathbf{S}}_{K,\boldsymbol{n}}-\sup_{t\in\mathbb{R}}
\mathbb{B}_{K,\boldsymbol{n}}^*(t)\right|=O\big(\phi_K^*(\boldsymbol{n})\big)
\quad\text{and}\quad\left|\overline{\mathbf{T}}_{K,\boldsymbol{n}}-\int_{\mathbb{R}}
\mathbb{B}_{K,\boldsymbol{n}}^*(t)\,dF_0(t)\right|=O\big(\phi_K^*(\boldsymbol{n})\big).
$$
\end{corollary}

\section{The change-point problem}\label{sectionchnae}

Here and elsewhere, $\lfloor t\rfloor$ denotes the largest integer not
exceeding $t$. In many practical applications, we assume the structural
stability of statistical models and this fundamental assumption needs to be
tested before it can be applied. This is called the analysis of structural
breaks, or change-points, which has led to the development of a variety of
theoretical and practical results. For good sources of references to research
literature in this area along with statistical applications, the reader may
consult \cite{Brodsky1993}, \cite{Csorgo1997} and \cite{ChenGupta2000}.
For recent references on the subject we may refer, among many others, to
\cite{Bouzebda2012MMS}, \cite{Alexander2012}, \cite{Julian2013}, \cite{Lajos2014},
\cite{Bouzebda2014AN} and \cite{Bouzebda2014MMS}.

In this section, we deal with testing changes in d.f.'s for a sequence of
independent real-valued r.v.'s $X_1,\ldots,X_n$. The corresponding null
hypothesis that we want to test is
$$
\mathcal{H}_0'': X_1,\ldots,X_n \text{~~have d.f.~~}F.
$$
As frequently done, the behavior of the derived tests will be investigated
under the alternative hypothesis of a single change-point
$$
\mathcal{H}_1'': \exists ~k^*\in\{1,\ldots,n-1\}\text{~~such that~~}
X_1,\ldots, X_{k^*}\text{~~have d.f.~~}F \text{~~and~~}X_{k^*+1},\ldots,X_n\text{~~have d.f.~~}G.
$$
The d.f.'s $F$ and $G$ are assumed to be continuous. The critical integer
$k^*$ can be written as $\lfloor ns\rfloor$ for a certain $s\in[0,1)$. Then,
testing the null hypothesis $\mathcal{H}_0''$ can be based on functionals of
the following process: set, for each $n\in\mathbb{N}^*$,
\begin{equation}\label{alphatilde}
\widetilde{\alpha}_n(s,t):=\frac{\lfloor ns\rfloor(n
-\lfloor ns\rfloor)}{n^{3/2}}\left(\overline{\mathbb{F}}_{\lfloor ns\rfloor}^-(t)
-\overline{\mathbb{F}}_{n-\lfloor ns\rfloor}^+(t)\right)
\quad\text{for}\quad s\in[0,1],\,t\in \mathbb{R},
\end{equation}
where $\overline{\mathbb{F}}_k^-$ is the integrated empirical d.f. based upon
the $k$ first observations and $\overline{\mathbb{F}}_{n-k}^+$ is that based upon
the $(n-k)$ last ones. In (\ref{alphatilde}) we extend the definition of
$\overline{\mathbb{F}}_k^-$ and $\overline{\mathbb{F}}_k^+$ to the case where $k=0$ by
setting $\overline{\mathbb{F}}_0^-=\overline{\mathbb{F}}_0^+=0$, so that
$\widetilde{\alpha}_n(s,t)=0$ if $s\in(0,1/n)$.

We can define the r.v.'s $X_1,\ldots ,X_{\lfloor ns\rfloor}$ and
$X_{\lfloor ns\rfloor+1},\ldots ,X_n$ on a probability space on which we can
simultaneously construct two Kiefer processes
$\{\mathbb{K}_1(s,u): s\in\mathbb{R}, u\in[0,1]\}$ and
$\{\mathbb{K}_2(s,u): s\in\mathbb{R}, u\in[0,1]\}$ such that the ``restricted''
processes $\{\mathbb{K}_1(s,u): s\in[1,n/2], u\in[0,1]\}$ and
$\{\mathbb{K}_2(s,u): s\in[n/2,n], u\in[0,1]\}$ are independent. It turns out
that a natural approximation of
$\big\{\widetilde{\alpha}_n:n\in\mathbb{N}^*\big\}$
is given by the sequence of Gaussian processes $\big\{\overline{\mathbb{K}}_n(s,F(t)) :
s\in[0,1], t\in\mathbb{R}, n\in\mathbb{N}^*\big\}$ defined by
$$
\overline{\mathbb{K}}_n(s,u):=\!
\begin{cases}
\frac{u}{\sqrt{n}}\big[\,\mathbb{K}_2(\lfloor ns \rfloor,u)
-s(\mathbb{K}_1(\lfloor n/2\rfloor,u)
+\mathbb{K}_2(\lfloor n/2\rfloor,u))\big]
& \text{for } s\in\!\big[0,\frac{1}{2}\big], u\in[0,1] \\[1ex]
\frac{u}{\sqrt{n}}\big[\!-\mathbb{K}_1(\lfloor n(1-s) \rfloor,u)
+(1-s)(\mathbb{K}_1(\lfloor n/2\rfloor,u)
+\mathbb{K}_2(\lfloor n/2\rfloor,u))\big]
& \text{for } s\in\!\big[\frac{1}{2},1\big], u\in[0,1].
\end{cases}
$$
More precisely, we have the following result.
\begin{theorem}\label{theorem1}
On a suitable probability space, it is possible to define
$\big\{\widetilde{\alpha}_n:n\in\mathbb{N}^*\big\}$,
jointly with a sequence of Gaussian processes
$\{\overline{\mathbb{K}}_n:n\in\mathbb{N}^*\}$ as above,
such that, under $\mathcal{H}_0''$, with probability $1$, as $n \to\infty$,
$$
\sup_{s\in[0,1]}\sup_{t\in \mathbb{R}}\left| \widetilde{\alpha}_n(s,t)
-\overline{\mathbb{K}}_n(s,F(t)) \right|=O\!\left(\frac{(\log n)^2}{\sqrt{n}}\right)\!.
$$
\end{theorem}
According to \cite{Cosorgo1997}, a way to test change-point is to use the
following statistics:
\begin{equation}\label{proempri}
\tau_n:=\sup_{s\in[0,1]}\sup_{t\in \mathbb{R}}
\left|\widetilde{\alpha}_n(s,t)\right|.
\end{equation}

The corollary below is a consequence of Theorem~\ref{theorem1}.
\begin{corollary}\label{coroltau}
If $\mathcal{H}_0''$ holds true, then we have the convergence in distribution,
as $n\to\infty$,
$$
\tau_n\stackrel{\mathcal{L}}{\longrightarrow} \sup_{s,u\in[0,1]}
\left|\overline{\mathbb{K}}(s,u)\right|,
$$
where $\overline{\mathbb{K}}=\big\{\overline{\mathbb{K}}(s,u): s,u\in[0,1]\big\}$
is a Gaussian process with mean zero and covariance function
$$
\mathbb{E}\left( \overline{\mathbb{K}}(s,u)\overline{\mathbb{K}}(s',u')\right)
=uu'(u\wedge u'-uu')(s\wedge s'-ss').
$$
\end{corollary}
One has $\overline{\mathbb{K}}(s,u)
= u\,\overset{\text{\tiny o}}{\mathbb{K}}(s,u)$ where
$\overset{\text{\tiny o}}{\mathbb{K}}$ is a tied-down Kiefer process.
We refer to \cite{Csorgo1997} for more details on the process $\overline{\mathbb{K}}$.

Actually, according to \cite{Cosorgo1997}, the most appropriate way to test
change-point is to use the following weighted statistic:
\begin{equation}\label{proempri2}
\tau_{n,w}:=\sup_{s\in[0,1]}\sup_{t\in\mathbb{R}}
\frac{\big|\widetilde{\alpha}_n(s,t)\big|}{w\left(\lfloor ns\rfloor /n\right)}
\end{equation}
where $w$ is a positive function defined on $(0,1)$, increasing in a neighborhood
of zero and decreasing in a neighborhood of one, satisfying the condition
$$
I(w,\boldsymbol{\varepsilon}):= \int_0^1\exp\left(-\frac{\boldsymbol{\varepsilon}
w^2(s)}{s(1-s)}\right)\,\frac{ds}{s(1-s)} <\infty
$$
for some constant $\boldsymbol{\varepsilon}>0$. For a history and further
applications of $I(w,\boldsymbol{\varepsilon})$, we refer to \cite{Csorgoho1993},
Chapter 4. From \cite{Szyszkowicz(1992)}, an example of such function $w$ is
given by
$$
w(t):=\left( t(1-t)\log \log \frac{1}{t(1-t)}\right)^{\!1/2}
\quad\text{for}\quad t\in(0,1).
$$
By using similar techniques to those which are developed in \cite{Csorgo1997},
one may show that, as $n\to\infty$,
$$
\tau_{n,w}\stackrel{\mathcal{L}}{\longrightarrow} \sup_{s,u\in[0,1]}
\frac{\big|\overline{\mathbb{K}}(s,u)\big|}{w(s)}.
$$
For more details, we refer to \cite{Bouzebda2014AN}.
\begin{remark}
As in \cite{Szysz(1994)}, we mention that the statistic given by
(\ref{proempri}) should be more powerful for detecting changes that occur in
the middle, i.e., near $n/2$, where $k/n(1-k/n)$ reaches its maximum, than for
the ones occurring near the end points. The advantage of using the weighted
statistic defined in (\ref{proempri2}) is the detection of changes that occur
near the end points, while retaining the sensitivity to possible changes in
the middle as well.
\end{remark}

We hope that the results presented in Sections~\ref{sectiontwo} and
\ref{sectionchnae} will be the prototypes of other various applications.

\section{Strong approximation of the integrated empirical process when
parameters are estimated}\label{Section5}

In this section, we are interested in the strong approximation of the
integrated empirical process when parameters are estimated. Our approach is in
the same spirit of \cite{Burke-Csorgo}. Let us introduce, for each
$n\in\mathbb{N}^*$, the \emph{integrated estimated empirical process}
$\overline{\widehat{\alpha}}_n$:
\begin{equation}\label{alphabarchap}
\overline{\widehat{\alpha}}_n(t):=\sqrt{n}\left( \overline{\mathbb{F}}_n(t)
-\overline{F}\big(t,\widehat{\boldsymbol{\theta}}_n\big)\right)
\quad\text{for}\quad t\in\mathbb{R},
\end{equation}
where $\big\{\widehat{\boldsymbol{\theta}}_n: n\in\mathbb{N}^* \big\}$ is a
sequence of estimators of a parameter $\boldsymbol{\theta}$ from a family of
d.f.'s $\{ F(t,\boldsymbol{\theta}): t\in\mathbb{R}, \,\boldsymbol{\theta}
\in\boldsymbol{\Theta}\}$ ($\boldsymbol{\Theta}$ being a subset of
$\mathbb{R}^d$ and $d$ a fixed positive integer) related to a sequence of
i.i.d. r.v.'s $\{X_i:i\in\mathbb{N}^*\}$. Let us mention that a general study
of the weak convergence of the estimated empirical process was carried out by
\cite{Durbin1973}. For a more recent reference, we may refer to \cite{Genz2006}
where the authors investigated the empirical processes with estimated
parameters under auxiliary information and provided some results regarding the
bootstrap in order to evaluate the limiting laws.

Let us introduce some notations.
\begin{enumerate}[label=(\thesection.\arabic*)]
\item
The transpose of a vector $V$ of $\mathbb{R}^d$ will be denoted by $V^\top$.
\item
The norm $\| \cdot\|$ on $\mathbb{R}^d$ is defined by
$$
\| (y_1,\ldots,y_d)\|:=\max_{1\leq i \leq d}|y_i|.
$$
\item
For a function $(t,\boldsymbol{\theta})\mapsto g(t,\boldsymbol{\theta})$ where
$\boldsymbol{\theta}=(\theta_1, \ldots,\theta_d)\in \mathbb{R}^d$,
$\nabla_{\boldsymbol{\theta}}g(t,\boldsymbol{\theta}_0)$
denotes the vector in $\mathbb{R}^d$ of partial derivatives
$\big((\partial g/\partial\theta_1) (t,\boldsymbol{\theta}),
\ldots,(\partial g/\partial\theta_d) (t,\boldsymbol{\theta}))\big)$
evaluated at $\boldsymbol{\theta}=\boldsymbol{\theta}_0$,
and $\nabla_{\boldsymbol{\theta}}^2g(t,\boldsymbol{\theta})$
denotes the $d\times d$ matrix of second order partial derivatives
$\big((\partial^2g/\partial\theta_i\partial\theta_j) (t,\boldsymbol{\theta}))
\big)_{1\leq i,j\leq d}$.
\item
For a vector $V=(v_1,\dots,v_d)\in\mathbb{R}^d$, $\int V$ denotes
the vector $\left(\int v_1,\dots,\int v_d\right)$.
\end{enumerate}

Next, we write out the set of all conditions (those of \cite{Burke-Csorgo})
which will be used in the sequel.
\begin{enumerate}[label=(\roman*)]
\item
The estimator $\widehat{\boldsymbol{\theta}}_n$ admits the following form:
for each $n\in\mathbb{N}^*$,
$$
\sqrt{n}\left(\widehat{\boldsymbol{\theta}}_n-\boldsymbol{\theta}_0\right)
=\frac{1}{\sqrt{n}}\sum_{i=1}^n l(X_i,\boldsymbol{\theta}_0)
+\boldsymbol{\varepsilon}_n,
$$
where $\boldsymbol{\theta}_0$ is the theoretical true value of
$\boldsymbol{\theta}$, $l(\cdot,\boldsymbol{\theta}_0)$ is a measurable
$d$-dimensional vector-valued function, and $\boldsymbol{\varepsilon}_n$
converges to zero as $n\to\infty$ in a manner to be specified later on.
Notice that
$$
\frac{1}{\sqrt{n}}\sum_{i=1}^n l(X_i,\boldsymbol{\theta}_0)
=\sqrt{n}\int_{-\infty}^t l(s,\boldsymbol{\theta}_0)\,d\mathbb{F}_n(s).
$$

\item
The mean value of $l(X_i,\boldsymbol{\theta}_0)$ vanishes:
$\mathbb{E}\!\left(l(X_i,\boldsymbol{\theta}_0)\right)=0$.
\item
The matrix $M(\boldsymbol{\theta}_0):=\mathbb{E}\!
\left(l(X_i,\boldsymbol{\theta}_0)^\top l(X_i,\boldsymbol{\theta}_0)\right)$
is a finite nonnegative definite $d\times d$ matrix.
\item
The vector-valued function $(t,\boldsymbol{\theta})\mapsto
\nabla_{\boldsymbol{\theta}}F(t,\boldsymbol{\theta})$ is uniformly
continuous in $t\in\mathbb{R}$ and $\boldsymbol{\theta} \in\mathbf{V}$, where
$\mathbf{V}$ is the closure of a given neighborhood of $\boldsymbol{\theta}_0$.
\item
Each component of the vector-valued function $t\mapsto
l(t,\boldsymbol{\theta}_0)$ is of bounded variation in $t$ on each finite interval of $\mathbb{R}$.
\item
The vector-valued function $t\mapsto\nabla_{\boldsymbol{\theta}}
F(t,\boldsymbol{\theta}_0)$ is uniformly bounded in
$t\in\mathbb{R}$, and the vector-valued function
$(t,\boldsymbol{\theta})\mapsto \nabla^2_{\boldsymbol{\theta}}F(t,\boldsymbol{\theta})$
is uniformly bounded in $t\in\mathbb{R}$ and $\boldsymbol{\theta} \in\mathbf{V}$.
\item
Set
$$
\ell(s,\boldsymbol{\theta}_0):=l\!\left(F^{-1}(s,\boldsymbol{\theta}_0),\boldsymbol{\theta}_0\right)
\quad\text{for}\quad s\in(0,1)
$$
where
$$
F^{-1}(s,\boldsymbol{\theta}_0)=\inf\{ t\in\mathbb{R}:F(t,\boldsymbol{\theta}_0)\geq s\}.
$$
The limiting relations below hold:
$$
\lim_{s\searrow 0} \sqrt{s \log \log (1/s)}\,
\left\|\ell(s,\boldsymbol{\theta}_0)\right\| =0
\quad\text{and}\quad\lim_{s\nearrow 1} \sqrt{(1-s) \log \log [1/(1-s)]}\,
\left\|\ell(s,\boldsymbol{\theta}_0)\right\| =0,
$$
\item
Set
$$
\ell'_s(s,\boldsymbol{\theta}_0):=\frac{\partial\ell}{\partial s}(s,\boldsymbol{\theta}_0)
\quad\text{for}\quad s\in(0,1).
$$
The partial derivative $\ell'_s(s,\boldsymbol{\theta}_0)$ exist for every $s\in (0,1)$ and
the bounds below hold: there is a positive constant $C$ such that
$$
s\left\|\ell'_s(s,\boldsymbol{\theta}_0)\right\| \leq C \quad\text{for all }
s\in\big(0,\textstyle{\frac{1}{2}}\big)\quad\text{and}\quad
(1-s)\left\|\ell'_s(s,\boldsymbol{\theta}_0)\right\| \leq C \quad\text{for all }
s\in\big(\textstyle{\frac{1}{2}},1\big).
$$
\end{enumerate}

Now, we state an analogous result to Theorem 3.1 of \cite{Burke-Csorgo}. For each
$n\in\mathbb{N}^*$, let $G_n=\{G_n(t):t\in\mathbb{R}\}$ be the process defined by
\begin{align*}
G_n(t)
&
:=\frac{1}{\sqrt{n}}\left( \mathbb{K}(n,F(t,\boldsymbol{\theta}_0))
- \left(\int_{\mathbb{R}} l(s,\boldsymbol{\theta}_0)\,
d_s\mathbb{K}(n,F(s,\boldsymbol{\theta}_0))\right)
\nabla_{\boldsymbol{\theta}}F(t,\boldsymbol{\theta}_0)^\top\right)
\\
&
=\frac{1}{\sqrt{n}}\left(\mathbb{K}(n,F(t,\boldsymbol{\theta}_0))
- \mathbf{W}(n)\nabla_{\boldsymbol{\theta}}F(t,\boldsymbol{\theta}_0)^\top\right)
\quad\text{for}\quad t\in\mathbb{R},
\end{align*}
where we set
$$
\mathbf{W}(\tau):=\int_{\mathbb{R}} l(s,\boldsymbol{\theta}_0)
\,d_s\mathbb{K}(\tau,F(s,\boldsymbol{\theta}_0))
\quad\text{for}\quad \tau\geq 0.
$$
The process $\{\mathbf{W}(\tau): \tau\geq 0\}$ is a $d$-dimensional Brownian
motion with a covariance matrix of rank that of $M(\boldsymbol{\theta}_0)$.
The estimated empirical process given by $\overline{\widehat{\alpha}}_n$ defined
by (\ref{alphabarchap}) will be approximated by the sequence of processes
$\overline{G}_n:=\{F(t,\boldsymbol{\theta}_0)\,G_n(t):t\in\mathbb{R}\}$. Set
$$
\overline{\boldsymbol{\varepsilon}}_n:=\sup_{t\in \mathbb{R}}
\big| \overline{\widehat{\alpha}}_n(t)-\overline{G}_n(t)\big|\!.
$$

\begin{theorem}\label{theoremBurkeetal}
Suppose that the sequence of estimators $\big\{\widehat{\boldsymbol{\theta}}_n:
n\in\mathbb{N}^*\big\}$ satisfies Conditions {\rm (i), (ii)} and {\rm (iii)}.
Then, as $n\to\infty$,
\begin{enumerate}[label=(\alph*)]
\item
$\overline{\boldsymbol{\varepsilon}}_n\stackrel{\mathbb{P}}{\longrightarrow} 0$
if Conditions {\rm (iv), (v)} hold and $\boldsymbol{\varepsilon}_n
\stackrel{\mathbb{P}}{\longrightarrow} 0$;
\item
$\overline{\boldsymbol{\varepsilon}}_n\stackrel{\text{a.s.}}{\longrightarrow} 0$
if Conditions {\rm (vi)--(viii)} hold and $\boldsymbol{\varepsilon}_n
\stackrel{\text{a.s.}}{\longrightarrow} 0$;
\item
$\overline{\boldsymbol{\varepsilon}}_n=O(\max(h(n),n^{-\epsilon}))$
for some $\epsilon>0$ if Conditions {\rm(vi)--(viii)} hold and
$\boldsymbol{\varepsilon}_n=O(h(n))$ for some function $h$ satisfying $h(n)>0$
and $h(n)\to 0.$
\end{enumerate}
\end{theorem}
The limiting Gaussian process $\overline{G}_n$ of Theorem~\ref{theoremBurkeetal}
depends crucially on $F$ and also on the true theoretical value
$\boldsymbol{\theta}_0$. In general, Theorem~\ref{theoremBurkeetal} cannot be
used to test the composite hypothesis :
$$
F\in\{ F(t,\boldsymbol{\theta}): t\in\mathbb{R}, \boldsymbol{\theta}\in\boldsymbol{\Theta}\}.
$$
In order to circumvent this problem, \cite{Burke-Csorgo} proposed an
approximate solution, they introduce another process:
$$
\widehat{G}_n(t):=\frac{1}{\sqrt{n}}\left( \mathbb{K}\big(n,F\big(t,
\widehat{\boldsymbol{\theta}}_n\big)\big)
- \mathbf{W}(n)\nabla_{\boldsymbol{\theta}}F\big(t,
\widehat{\boldsymbol{\theta}}_n\big)^\top\right)\!.
$$
Under some regularity conditions, \cite{Burke-Csorgo} show that (see Theorem~3.2
therein), as $n\to\infty,$
$$
\sup_{t\in\mathbb{R}}\left|\widehat{G}_n(t)-G_n(t)\right|
\stackrel{\mathbb{P}}{\longrightarrow }0.
$$
Setting $\overline{\widehat{G}}_n(t) = F(t,\widehat{\boldsymbol{\theta}}_n)\,
\widehat{G}_n(t)$, one can show that,
as $n\to\infty$,
\begin{equation}\label{eqreferrrz11}
\sup_{t\in\mathbb{R}}\left|\overline{\widehat{G}}_n(t)-
\overline{G}_n(t)\right|\stackrel{\mathbb{P}}{\longrightarrow }0.
\end{equation}
Consequently, we have, as $n\to\infty$,
$$
\sup_{t\in \mathbb{R}}\left| \overline{\widehat{\alpha}}_n(t)-\overline{\widehat{G}}_n(t)\right
|\stackrel{\mathbb{P}}{\longrightarrow }0.
$$

\section{Local time of the integrated empirical process}\label{section343}

In this section, we are mainly concerned with the behavior of the local time of
the integrated empirical process. This characterization is possible by making use
of a representation provided by \cite{HenzeNikitin2002} that expresses the
integrated empirical process in terms of a partial sums process, see (\ref{sj})
below. Recall the definition of the process given in (\ref{beta2}). Set
$$
\overline{\beta}_n(u):=\int_0^u\beta_n(v)\,dv.
$$
In this part, we focus on the particular r.v. $\overline{\mathcal{A}}_n
=\overline{\beta}_n(1)$. According to \cite{HenzeNikitin2002} p.~185, we have
the following representation for $\overline{\mathcal{A}}_n$:
$$
\overline{\mathcal{A}}_n=\sqrt{n}\left(\frac{1}{n}\sum_{i=1}^n(1-U_i)-\frac{1}{2} \right)
=\frac{S_n}{\sqrt{n}}
$$
where $\{S_n:n\in\mathbb{N}^*\}$ is the following partial sums
process where the summands are i.i.d. r.v.'s with mean zero and values
in the interval $I:=[-1/2,1/2]$:
\begin{equation}\label{sj}
S_n:=\sum_{i=1}^n\left(\frac{1}{2}-U_i \right)\!.
\end{equation}
Notice that we are dealing with a sum of strongly non-lattice r.v.'s as, i.e.,
in p.~210 of \cite{BassKhoshnevisan1993}. Indeed, we easily check that the
characteristic function $\chi$ of the $(1/2-U_i)$'s,
namely
$$
\chi(z):=\frac{\sin(z/2)}{z/2},
$$
satisfies the conditions
$$
\forall z\in\mathbb{R}^*,\,|\chi(z)|< 1\quad\text{and}\quad \limsup_{\left|z\right|\to\infty}|\chi(z)|< 1.
$$

Next, we define the local time
\begin{equation}\label{localtime}
\lambda(x,n):=\sum_{i=1}^n\mathbbm{1}_I(S_i-x)\quad\text{for}
\quad x\in \mathbb{R},\, n\in\mathbb{N}^*.
\end{equation}
The local time $\lambda(x,n)$ represents the number of visits of the
random walk $\{S_n:n\in\mathbb{N}^*\}$ in the neighborhood
$x+I$ of $x$ up to discrete time $n$. \cite{Revesz1981} proved that if the
random walk is symmetric, then, on some enlarged probability space supporting
also a Brownian motion with local time $l(x,t)$ defined as in (\ref{Lb}) below,
we have, for any $\varepsilon>0$, almost surely, as $n\to\infty$,
$$
\sup_{x\in\mathbb{Z}}|\lambda(x,n)-l(x,n)|=O\!\left(n^{1/4+\varepsilon}\right)\!.
$$
Actually, \cite{BassKhoshnevisan1993} obtained the more precise estimate (see Theorem 4.5
therein) below.
\begin{theo}
We have, with probability~$1$, as $n \to\infty$,
\begin{equation}\label{BaKh1}
\sup_{x\in\mathbb{R}}\sup_{t\in[0,1]}|\lambda(x,\lfloor nt\rfloor)-l(x,nt)|
=O\!\left(n^{1/4}(\log n)^{1/2}(\log \log n)^{1/4}\right)\!.
\end{equation}
\end{theo}

In the present context, our aim is to obtain the rate of the approximation
of the self-intersection local time
$$
L_n(t):=\sum_{1\leq i<j \leq \lfloor nt\rfloor }\int_{\mathbb{R}}
\mathbbm{1}_I(S_i-x)\mathbbm{1}_I(S_j-x)\,dx
$$
by the integrated local time of some standard Wiener process.
The quantity $L_n(t)$ enumerates in a certain manner the couples $(i,j)$
of distinct and ordered indices up to time $\lfloor nt\rfloor$ such that
$S_i-S_j$ is less than the diameter of $I$.

To this aim, we recall that, if $\{\mathbb{W}(t): t\geq 0\}$ is the standard Wiener
process with $\mathbb{W}(0)=0$, then its local time process
$\{l(x,t): t \geq 0, x\in\mathbb{R}\}$ is defined as
\begin{equation}\label{Lb}
l(x,t):=\lim_{\boldsymbol{\varepsilon} \downarrow 0}\frac{1}{2\boldsymbol{\varepsilon}}
\int_0^t\mathbbm{1}_{\left\{ x-\boldsymbol{\varepsilon} \leq \mathbb{W}(s) \leq x
+\boldsymbol{\varepsilon} \right\}}\,ds\quad\text{for}\quad x\in\mathbb{R},\,t \geq 0.
\end{equation}
Let us introduce the corresponding normalized local time $l_n$:
$$
l_n(x,t):=\frac{1}{\sqrt{n}}\,l\!\left(\sqrt{n}\,x,\lfloor nt\rfloor\right)
\quad\text{for}\quad x\in \mathbb{R},\,t \geq 0,\, n\in\mathbb{N}^*.
$$
\begin{theorem}\label{Res1}
We have, with probability~$1$, as $n \to\infty$,
$$
\sup_{t\in[0,1]} \left|L_n(t)
-\frac{1}{2}\,n^{3/2}\int_{\mathbb{R}} l_n(x,t)^2\,dx\right|
=O\!\left(n^{5/4}(\log n)^{1/2}(\log \log n)^{1/4}\right)\!.
$$
\end{theorem}

Lemma~2.2 of \cite{HuKhoshnevisan2010} stipulates, almost surely, as $n \to\infty$,
$$
\int_{\mathbb{R}}l(x,n)^2\,dx=n^{3/2+o(1)}.
$$
Actually, when looking at the proof of this lemma, one can prove
more precisely that, almost surely, there exist two positive constants
$\kappa_1'$ and $\kappa_2'$ such that, as $n \to\infty$,
\begin{equation}\label{inequalityln}
\kappa_1'\,\frac{n^{3/2}}{\sqrt{\log\log n}}\leq\int_{\mathbb{R}}l(x,n)^2\,dx
\leq \kappa_2'\,n^{3/2}\sqrt{\log\log n}.
\end{equation}
we obtain the following result.
\begin{corollary}\label{corLn}
We have, with probability $1$, for any $t\in(0,1]$, there exist two positive
constants $\kappa_1$ and $\kappa_2$ such that, almost surely, for large enough $n$,
\begin{equation}\label{inequalityLn}
\kappa_1\,\frac{\lfloor nt\rfloor^{3/2}}{\sqrt{\log\log n}}\leq L_n(t)
\leq \kappa_2\,\lfloor nt\rfloor^{3/2}\sqrt{\log\log n}.
\end{equation}
In particular, for any $t\in(0,1]$, almost surely, as $n \to\infty$,
$$
L_n(t)=\frac{1}{2}\,\lfloor nt\rfloor^{3/2+o(1)}.
$$
\end{corollary}

\section{Mathematical developments}\label{appe}

This section is devoted to the detailed proofs of our results. The previously
displayed notations continue to be used in the sequel.
\\

\subsection{Proof of Formula (\ref{beta21})}

For the completeness of the paper, we provide here a proof of
(\ref{beta21}) inspired by \cite{HenzeNikitin2002}. The function $\overline{\mathbb{F}}_n$ can
be computed as follows: for any $n\in\mathbb{N}^*$ and any $t\in\mathbb{R}$,
$$
\overline{\mathbb{F}}_n(t)=\frac{1}{n}\sum_{j=1}^n\int_{-\infty}^t\mathbb{F}_n(x)\,
d\mathbbm{1}_{\{X_j\leq x\}}=\sum_{j=1}^n\mathbb{F}_n(X_j)\mathbbm{1}_{\{X_j\leq t\}}
=\frac{1}{n^2}\sum_{1\leq i,j\leq n }\mathbbm{1}_{\{X_i\leq X_j\leq t\}}.
$$
Because of the hypothesis that the d.f. $F$ is continuous, the sampled variables
$X_1, X_2,\ldots,X_n$ are almost surely all different. Then, we can define with probability $1$
the order statistics
$$
X_{(1)} < X_{(2)} < \cdots< X_{(n)}
$$
associated with $X_1, X_2,\ldots,X_n$. Notice that the inequality $X_{(i)} \leq X_{(j)}$ is equivalent
to $i \leq j$, and that the event $\{X_{(j)} \leq t\}$ is equal to
$\{n\mathbb{F}_n(t) \leq j\}$. Hence, we can write that, with probability $1$,
for any $n\in\mathbb{N}^*$ and any $t\in\mathbb{R}$,
\begin{align*}
\overline{\mathbb{F}}_n(t)&=\frac{1}{n^2}\sum_{1\leq i,j\leq n }
\mathbbm{1}_{\{X_i\leq X_j\leq t\}}=\frac{1}{n^2}\sum_{1\leq i,j\leq n }
\mathbbm{1}_{\{X_{(i)}\leq X_{(j)}\leq t\}}
\\
&=\frac{1}{n^2}\# \left\{(i,j)\in \mathbb{N}^2: 1\leq i \leq j\leq n\mathbb{F}_n(t)\right\}\!.
\end{align*}
Now, we see that the foregoing cardinality is nothing but the number of combinations with repetitions
of two integers lying between $1$ and $n\mathbb{F}_n(t)$, which coincides with the number of
combinations without repetitions of two integers lying between $1$ and $n\mathbb{F}_n(t) + 1$, namely
$$
{{n\mathbb{F}_n(t)+1}\choose{2}}.
$$
This completes the proof of (\ref{beta21}). \hfill$\Box$\\

\subsection{Proof of Theorem~\ref{lem2}}

Making use of (\ref{beta21}) together with the elementary identity
$a^2 - b^2 = 2b(a-b) + (a -b)^2$, we can write, a.s.,
for any $t\in\mathbb{R},\,n\in\mathbb{N}^*$,
\begin{align}\label{H-N}
\nonumber\overline{\alpha}_n(t)&=\frac{1}{2}\,\sqrt{n}\left(\mathbb{F}_n(t)^2-F(t)^2\right)
+\frac{1}{2\sqrt{n}}\,\mathbb{F}_n(t)
\\
\nonumber&=\sqrt{n}\,F(t)\left(\mathbb{F}_n(t)-F(t)\right)+\frac{1}{2}\,\sqrt{n}
\left(\mathbb{F}_n(t)-F(t)\right)^2+ \frac{1}{2\sqrt{n}}\,\mathbb{F}_n(t)
\\
&=F(t)\,\alpha_n(t)+\frac{1}{2\sqrt{n}}\,\alpha_n(t)^2+\frac{1}{2\sqrt{n}}\,\mathbb{F}_n(t).
\end{align}
Since $0 \leq \mathbb{F}_n(t) \leq 1$, it is clear that, a.s.,
for any $t\in\mathbb{R},\,n\in\mathbb{N}^*$,
$$
F(t)\,\alpha_n(t)+\frac{1}{2\sqrt{n}}\,\alpha_n(t)^2\leq \overline{\alpha}_n(t)
\leq F(t)\,\alpha_n(t)+\frac{1}{2\sqrt{n}}\,\alpha_n(t)^2+ \frac{1}{2\sqrt{n}}.
$$
Consequently, since $0 \leq F(t) \leq 1$, we have (by roughly bounding $1/2$ by $1$),
a.s., for any $n\in\mathbb{N}^*$,
$$
\sup_{t\in \mathbb{R}}\big|\overline{\alpha}_n(t)-F(t)\,\mathbb{B}_n(F(t))\big|\leq
\sup_{t\in \mathbb{R}}\big|\alpha_n(t)-\mathbb{B}_n(F(t))\big|
+\frac{1}{\sqrt{n}}\sup_{t\in \mathbb{R}}\alpha_n(t)^2+ \frac{1}{\sqrt{n}},
$$
from which we deduce
\begin{align*}
\lqn{\mathbb{P}\!\left\{ \sup_{t\in \mathbb{R}}\big|\overline{\alpha}_n(t)
-F(t)\,\mathbb{B}_n(F(t))\big| > \frac{1}{\sqrt{n}}\,(A\log n+x)\right\}}
&
\leq \mathbb{P}\!\left\{ \sup_{t\in \mathbb{R}}\big|\alpha_n(t)-\mathbb{B}_n(F(t))\big|
+\frac{1}{\sqrt{n}}\sup_{t\in \mathbb{R}}\alpha_n(t)^2 > \frac{1}{\sqrt{n}}\,(A\log n+x-1)\right\}\!.
\end{align*}
Now, using the elementary inequality $\mathbb{P}\{X+Y>a\} \leq \mathbb{P}\{X>b\}+\mathbb{P}\{Y>a-b\}$
which is valid for any r.v.'s $X, Y$ and any real numbers $a, b$, we obtain
\begin{align}
\lqn{\mathbb{P}\!\left\{ \sup_{t\in \mathbb{R}} \big|\overline{\alpha}_n(t)
-F(t)\,\mathbb{B}_n(F(t))\big| > \frac{1}{\sqrt{n}}\,(A\log n+x)\right\}}
&
\leq \mathbb{P}\!\left\{ \sup_{t\in \mathbb{R}}\big|\alpha_n(t)
-\mathbb{B}_n(F(t))\big|> \frac{1}{\sqrt{n}}\,(A\log n+\frac{x}{2}-1)\right\}+
\mathbb{P}\!\left\{ \sup_{t\in \mathbb{R}}\alpha_n(t)^2 > \frac{x}{2}\right\}\!.
\label{appp1ref}
\end{align}
The inequality of \cite{Dvoretzky1956} stipulates that there
exists a positive constant $c_4$ such that
\begin{equation}\label{equaDvoretzky1956}
\mathbb{P}\!\left\{\sup_{t\in \mathbb{R}}
\left|\mathbb{F}_n(t)-F(t)\right|>\sqrt{\frac{x}{n}}\,\right\}\leq c_4\,e^{-2x}
\end{equation}
for all $x>0$ and all $n\in\mathbb{N}^*$. An application of (\ref{equaDvoretzky1956})
gives for any $x>1$ and all $n\in\mathbb{N}^*$
\begin{equation}\label{appp1}
\mathbb{P}\!\left\{\sup_{t\in \mathbb{R}}\alpha_n(t)^2 >\frac{x}{2}\right\} \leq c_4\,e^{-x}.
\end{equation}
Now, by putting (\ref{major}) and (\ref{appp1}) into (\ref{appp1ref}), we complete the proof of
Theorem~\ref{lem2}. \hfill$\Box$\\

\subsection{Proof of Corollary~\ref{corol22}}

The functional $\Phi$ being Lipschitz, there exists a positive constant $L$
such that, for any functions $u,v$,
$|\Phi(u)-\Phi(v)|\leq L \sup_{t\in\mathbb{R}} |u(t)-v(t)|$, inequality
that we will use in the form
\begin{equation}\label{lip}
\Phi(v)-L \sup_{t\in\mathbb{R}} |u(t)-v(t)|\leq \Phi(u) \leq \Phi(v)
+L \sup_{t\in\mathbb{R}} |u(t)-v(t)|.
\end{equation}
Let us choose for $u,v$ the processes $U_n:=\overline{\alpha}_n(\cdot)$ and
$V_n:=F(\cdot)\,\mathbb{B}_n(F(\cdot))$.
Applying the elementary inequality $|\mathbb{P}(A)-\mathbb{P}(B)|
\leq \mathbb{P}(A\backslash B) + \mathbb{P}(B\backslash A)$
to the events $A=\{\Phi(U_n)\leq x\}$ and $B=\{\Phi(V_n)\leq x\}$ provides,
for any $x\in\mathbb{R}$ and any $n\in\mathbb{N}^*$,
$$
\big|\mathbb{P} \{\Phi(U_n)\leq x\}-\mathbb{P} \{\Phi(V_n)\leq x\}\big|
\leq \mathbb{P} \{\Phi(U_n)\leq x\leq \Phi(V_n)\}
+\mathbb{P} \{\Phi(V_n)\leq x\leq \Phi(U_n)\}.
$$
By (\ref{lip}), we see that
\begin{align*}
\mathbb{P} \{\Phi(U_n)\leq x\leq \Phi(V_n)\}
&
\leq\mathbb{P} \!\left\{\Phi(V_n)-L \,\sup_{t\in\mathbb{R}} |U_n(t)-V_n(t)|\leq x\leq \Phi(V_n)\right\}\!,
\\
\mathbb{P} \{\Phi(V_n)\leq x\leq \Phi(U_n)\}
&
\leq\mathbb{P} \!\left\{\Phi(V_n)\leq x\leq \Phi(V_n)+L \,\sup_{t\in\mathbb{R}} |U_n(t)-V_n(t)|\right\}\!,
\end{align*}
from which we deduce, by addition, that
\begin{equation}\label{estimPhiUV}
\big|\mathbb{P} \{\Phi(U_n)\leq x\}-\mathbb{P} \{\Phi(V_n)\leq x\}\big|
\leq \mathbb{P} \!\left\{|\Phi(V_n)-x| \leq L \sup_{t\in\mathbb{R}} |U_n(t)-V_n(t)|\right\}\!.
\end{equation}

On the other hand, by choosing $x=c\,\log n$ for a suitable constant $c$ in
(\ref{estimation}) that will be specified below and putting $\epsilon_n:=(A+c)\log n/\sqrt{n}$,
we obtain the estimate below valid for large enough $n$:
$$
\mathbb{P}\!\left\{ \sup_{t\in\mathbb{R}} \left|U_n(t)-V_n(t)\right|
\geq \epsilon_n\right\}\leq \frac{B}{n^{cC}}.
$$
By choosing $c>1/(2C)$, we have
\begin{equation}\label{estimationbis}
\mathbb{P}\!\left\{ \sup_{t\in\mathbb{R}} \left|U_n(t)-V_n(t)\right|
\geq \epsilon_n\right\}=o\!\left(\frac{1}{\sqrt{n}}\right)\!.
\end{equation}
Now, by (\ref{estimPhiUV}), we write
\begin{align}
\big|\mathbb{P} \{\Phi(U_n)\leq x\}-\mathbb{P} \{\Phi(V_n)\leq x\}\big|=
&\;
\mathbb{P} \!\left\{ \sup_{t\in\mathbb{R}}
\left|U_n(t)-V_n(t)\right|<\epsilon_n, |\Phi(V_n)-x| \leq
L \sup_{t\in\mathbb{R}} |U_n(t)-V_n(t)|\right\}
\nonumber\\
&
+\mathbb{P}\!\left\{ \sup_{t\in\mathbb{R}} \left|U_n(t)-V_n(t)\right|
\geq \epsilon_n,|\Phi(V_n)-x| \leq
L \sup_{t\in\mathbb{R}} |U_n(t)-V_n(t)|\right\}
\nonumber\\
\leq&\;
\mathbb{P} \!\left\{|\Phi(V_n)-x| \leq L \epsilon_n\right\}
+\mathbb{P}\!\left\{ \sup_{t\in\mathbb{R}} \left|U_n(t)-V_n(t)\right|
\geq \epsilon_n\right\}\!.
\label{estimationter}
\end{align}
Noticing that the distribution of $\mathbb{B}_n$ does not depend on $n$,
which entails the equality
$\mathbb{P} \!\left\{|\Phi(V_n)-x| \leq L \epsilon_n\right\}
=\mathbb{P} \!\left\{|\Phi(V)-x| \leq L \epsilon_n\right\}$
where $V:=F(\cdot)\,\mathbb{B}(F(\cdot))$,
and recalling the assumption that the r.v.
$\Phi(V)$ admits a density function bounded by $M$ say,
we get that, for any $x\in\mathbb{R}$ and any $n\in\mathbb{N}^*$,
\begin{equation}\label{estimationquater}
\mathbb{P} \!\left\{|\Phi(V_n)-x| \leq L \epsilon_n\right\}\leq 2LM\epsilon_n.
\end{equation}
Finally, putting (\ref{estimationbis}) and  (\ref{estimationquater})
into (\ref{estimationter}) leads to (\ref{estimationfunctional}),
which completes the proof of Corollary~\ref{corol22}.
\hfill$\Box$

\subsection{Proof of Corollary~\ref{corollaryapprox}}

Applying (\ref{estimation}) to $x=c'\log n$ for some positive constant $c'$
that will be specified below yields, for large enough~$n$,
$$
\pi_n:=\mathbb{P}\!\left\{ \sup_{t\in\mathbb{R}} \big|\overline{\alpha}_n(t)
-F(t)\,\mathbb{B}_n(F(t))\big| \geq c_5\,\frac{\log n}{\sqrt{n}}\right\}
\leq \frac{B}{n^{c_6}}
$$
where we set $c_5=A+c'$ and $c_6=c'C$. Now, we choose $c'>1/C$ so that
the series $(\sum \pi_n)$ is convergent. Hence, by appealing to
Borel-Cantelli lemma, we get that
$$
\mathbb{P}\!\left(\limsup_{n\in\mathbb{N}^*}\left\{ \sup_{t\in\mathbb{R}}
\big|\overline{\alpha}_n(t)-F(t)\,\mathbb{B}_n(F(t))\big|
\geq c_5\,\frac{\log n}{\sqrt{n}}\right\}\right)=0
$$
which clearly implies Corollary~\ref{corollaryapprox}.
\hfill$\Box$

\subsection{Proof of Theorem~\ref{kieferapproximation}}

We write that, with probability $1$, for any $n\in\mathbb{N}^*$,
\begin{align}
\lqn{\max_{1\leq k\leq n}\sup_{t\in \mathbb{R}}\left|\sqrt{k}\,\overline{\alpha}_k(t)
-F(t)\,\mathbb{K}(k,F(t))\right|}
&
\leq \max_{1\leq k\leq n}\sup_{t\in \mathbb{R}}\left|\sqrt{k}\,F(t)\,\alpha_k(t)
+\frac{1}{2}\,\alpha_k^2(t)+\frac{1}{2}\,\mathbb{F}_k(t)-F(t)\,\mathbb{K}(k,F(t))\right|\nonumber
\nonumber\\
&\leq \max_{1\leq k\leq n}\sup_{t\in \mathbb{R}}\left|\sqrt{k}\,\alpha_k(t)
-\mathbb{K}(k,F(t))\right|+\max_{1\leq k\leq n}\sup_{t\in \mathbb{R}}\alpha_k^2(t)+1.
\label{referthe25}
\end{align}
In the last inequality above, we have used the facts that $0 \leq F(t) \leq 1$ and
$0 \leq \mathbb{F}_n(t) \leq 1$ for all $t\in \mathbb{R}$ and roughly bounded $1/2$ by $1$.
Now, by \cite{KMT1975}, we have, almost surely, as $n\to\infty$,
\begin{equation}\label{equation11KIE}
\max_{1\leq k\leq n}\sup_{t\in \mathbb{R}}\left|\sqrt{k}\,\alpha_k(t)-\mathbb{K}(k,F(t))\right|
=O\!\left((\log n)^2\right)\!.
\end{equation}
On the other hand, by using \cite{Chung1949}'s law of the iterated logarithm for the empirical process
which stipulates that
$$
\limsup_{n\to\infty} \frac{\sup_{t\in\mathbb{R}} |\alpha_n(t)|}
{\sqrt{\log \log n}}=\frac{1}{\sqrt{2}} \quad\text{a.s.},
$$
we see that, almost surely, as $n\to\infty$,
\begin{equation}\label{empililaeae}
\sup_{t\in\mathbb{R}} |\alpha_n(t)|=O\!\left(\!\sqrt{\log \log n}\,\right)\!.
\end{equation}
By putting (\ref{equation11KIE}) and (\ref{empililaeae}) into (\ref{referthe25}), we completes the proof
of Theorem~\ref{kieferapproximation}.\hfill$\Box$\\

\subsection{Proof of Corollary~\ref{ecooooo}}

We work under Hypothesis $\mathcal{H}_0$. Let us introduce the integrated
empirical process related to the d.f. $F_0$:
$$
\overline{\alpha}_{0,n}(t):=\sqrt{n}\left(\overline{\mathbb{F}}_n(t)-\overline{F}_0(t)\right)
\quad\text{for}\quad t\in\mathbb{R},\,n\in\mathbb{N}^*.
$$
By the triangular inequality, we plainly have
$$
\left|\overline{\mathbf{S}}_n-\sup_{t\in\mathbb{R}}
\big|F_0(t)\,\mathbb{B}_n(F_0(t))\big|\right|
\leq
\sup_{t\in\mathbb{R}}\big|\overline{\alpha}_{0,n}(t)
-F_0(t)\,\mathbb{B}_n(F_0(t))\big|
$$
from which together with (\ref{approx}) we deduce (\ref{kol}).
Similarly,
\begin{align}
\left|\overline{\mathbf{T}}_n-\int_{\mathbb{R}}
\big[F_0(t)\,\mathbb{B}_n(F_0(t))\big]^2\,dF_0(t)\right|
\leq&\;
\int_{\mathbb{R}}\left|\overline{\alpha}_{0,n}(t)^2
-\big[F_0(t)\,\mathbb{B}_n(F_0(t))\big]^2\right| dF_0(t)
\nonumber\\
\leq&\;
\sup_{t\in\mathbb{R}}\big|\overline{\alpha}_{0,n}(t)
-F_0(t)\,\mathbb{B}_n(F_0(t))\big|
\nonumber\\
&
\times \left(\sup_{t\in\mathbb{R}}\left|\overline{\alpha}_{0,n}(t)\right|
+\sup_{t\in\mathbb{R}}\big|F_0(t)\,\mathbb{B}_n(F_0(t))\big|\right)\!.
\label{vonmisesbis}
\end{align}
In the last inequality above appears the supremum
$$
\sup_{t\in\mathbb{R}}\big|F_0(t)\,\mathbb{B}_n(F_0(t))\big|
\leq\sup_{u\in[0,1]}\big|\mathbb{B}_n(u)\big|.
$$
On the other hand, by \cite{KMT1975}, on a suitable probability space, we can define
the empirical process $\{\beta_n:n\in\mathbb{N}^*\}$, in combination with the
sequence of Brownian bridges $\left\{\mathbb{B}_n:n\in\mathbb{N}^*\right\}$,
such that, with probability~$1$, as $n\to\infty$,
\begin{equation}\label{approxbis}
\sup_{u\in[0,1]} \left|\beta_n(u)- \mathbb{B}_n(u)\right|
=O\!\left(\frac{\log n}{\sqrt{n}}\right)\!.
\end{equation}
Actually, this is a consequence of (\ref{major}).
In view of (\ref{empililaeae}) and (\ref{approxbis}), one derives the following
bound: with probability~$1$,
as $n\to\infty$,
\begin{equation}\label{estimbridge}
\sup_{u\in[0,1]} |\mathbb{B}_n(u)| =O\!\left(\!\sqrt{\log \log n}\,\right)\!.
\end{equation}
Finally, by putting (\ref{approx}), (\ref{lli}) and (\ref{estimbridge}) into
(\ref{vonmisesbis}), we immediately deduce (\ref{vonmises}). The proof of
Corollary~\ref{ecooooo} is finished.
\hfill$\Box$

\subsection{Proof of Corollary~\ref{ximn}}

For each $m,n\in\mathbb{N}^*$, let $\overline{\alpha}_m^1$ and
$\overline{\alpha}_n^2$ denote the empirical
processes respectively associated with the samples $X_1,\ldots, X_m$ and $Y_1,\ldots,Y_n$.
By replacing $\overline{\mathbb{F}}_m(t)$ by $\overline{\alpha}_m^1(t)/\sqrt{m}+\overline{F}(t)$
and $\overline{\mathbb{G}}_n(t)$ by $\overline{\alpha}_n^2(t)/\sqrt{n}+\overline{G}(t)$,
using the binomial theorem and recalling that, under $\mathcal{H}_0'$, $F=G$,
we write
\begin{align*}
\overline{\boldsymbol{\xi}}_{m,n}^{(q)}(t)
&
=\sqrt{\frac{mn}{m+n}}\left[\left(\frac{\overline{\alpha}_m^1(t)}{\sqrt{m}}+\overline{F}(t)\right)^{\!\!q}
-\left(\frac{\overline{\alpha}_n^2(t)}{\sqrt{n}}+\overline{F}(t)\right)^{\!\!q}\right]
\\
&
=\sqrt{\frac{mn}{m+n}}\sum_{k=1}^q {q\choose k} \overline{F}(t)^{q-k}
\left[\left(\frac{\overline{\alpha}_m^1(t)}{\sqrt{m}}\right)^{\!\!k}
-\left(\frac{\overline{\alpha}_n^2(t)}{\sqrt{n}}\right)^{\!\!k}\right]
\\
&
=q\, \overline{F}(t)^{q-1} \left(\sqrt{\frac{n}{m+n}}\,\overline{\alpha}_m^1(t)
-\sqrt{\frac{m}{m+n}}\,\overline{\alpha}_n^2(t)\right) +\Delta_{m,n}(t)
\end{align*}
where we set
$$
\Delta_{m,n}(t)=\sqrt{\frac{mn}{m+n}}\sum_{k=2}^q {q\choose k} \overline{F}(t)^{q-k}
\left[\left(\frac{\overline{\alpha}_m^1(t)}{\sqrt{m}}\right)^{\!\!k}
-\left(\frac{\overline{\alpha}_n^2(t)}{\sqrt{n}}\right)^{\!\!k}\right]\!.
$$
By (\ref{lli}) and the fact that $0 \leq F(t) \leq 1$ for any $t\in\mathbb{R}$,
it is easily seen that, with probability $1$, as $m,n\to\infty$,
\begin{equation}\label{Delta}
\sup_{t\in\mathbb{R}}|\Delta_{m,n}(t)|=O\!\left(\frac{(\log\log m)^{q/2}}{\sqrt{m}}\right)
+O\!\left(\frac{(\log\log n)^{q/2}}{\sqrt{n}}\right)\!.
\end{equation}
On the other hand, by Corollary~\ref{corollaryapprox}, we can construct two sequences of
Brownian bridges $\big\{\mathbb{B}_m^1:m\in\mathbb{N}^*\big\}$
and $\big\{\mathbb{B}_n^2:n\in\mathbb{N}^*\big\}$ such that, with probability $1$, as $m,n\to\infty$,
\begin{equation}
\sup_{t\in\mathbb{R}} \left|\overline{\alpha}_m^1(t)-
F(t)\,\mathbb{B}_m^1(F(t))\right|=O\!\left(\frac{\log m}{\sqrt{m}}\right)\!,
\quad\sup_{t\in\mathbb{R}} \left|\overline{\alpha}_n^2(t)-
F(t)\,\mathbb{B}_n^2(F(t))\right|=O\!\left(\frac{\log n}{\sqrt{n}}\right)\!.
\label{ecart1}
\end{equation}
Setting $\mathbb{B}_{m,n}^{*(q)}$ as in Corollary~\ref{ximn}, we have
\begin{align}
\overline{\boldsymbol{\xi}}_{m,n}^{(q)}(t)-\mathbb{B}_{m,n}^{*(q)}(t)=
&\;
\frac{q}{2^{q-1}}\,F(t)^{2(q-1)}\left[
\sqrt{\frac{n}{m+n}}\left(\,\overline{\alpha}_m^1(t)-F(t)\,
\mathbb{B}_m^1(F(t))\right)\right.
\nonumber\\
&
\left.-\sqrt{\frac{m}{m+n}}\left(\,\overline{\alpha}_n^2(t)-F(t)\,
\mathbb{B}_n^2(F(t))\right)\right]+\Delta_{m,n}(t).
\label{ecart3}
\end{align}
By putting (\ref{Delta}) and (\ref{ecart1}) into (\ref{ecart3}),
we deduce the result announced in Corollary~\ref{ximn}.
\hfill$\Box$

\subsection{Proof of Theorem~\ref{theoremK}}

Let us introduce, for each $k\in\{1,\dots,K\}$, the integrated empirical
process associated with the d.f. $F^k$
$$
\overline{\alpha}_n^k(t):=\sqrt{n}\left(\overline{\mathbb{F}}_n^k(t)-\overline{F}^k(t)\right)
\quad\text{for}\quad t\in\mathbb{R},\,n\in\mathbb{N}^*.
$$
By recalling (\ref{DD}) and making use of the most well-known variance formula
$$
\sum_{k=1}^K n_k(x_k-\bar{x})^2
=\sum_{k=1}^K n_k(x_k-x_0)^2 -\boldsymbol{|n|}\left(\bar{x}-x_0\right)^2
$$
where we have denoted $\boldsymbol{|n|}=\sum_{k=1}^K n_k$
and $\bar{x}=\frac{1}{\boldsymbol{|n|}}\sum_{k=1}^K n_kx_k$
for $\boldsymbol{n}=(n_1,\dots,n_K)$, we rewrite
$\overline{\boldsymbol{\xi}}_{K,\boldsymbol{n}}(t)$ under Hypothesis
$\mathcal{H}_0^K$ as
\begin{align*}
\overline{\boldsymbol{\xi}}_{K,\boldsymbol{n}}(t)
&
=\sum_{k=1}^K n_k\!
\left(\overline{\mathbb{F}}_{n_k}^k(t)-\overline{F}_0(t)\right)^{\!2}
-\frac{1}{\boldsymbol{|n|}}\left(\sum_{k=1}^K n_k\!
\left(\overline{\mathbb{F}}_{n_k}^k(t)-\overline{F}_0(t)\right)\right)^{\!\!2}
\\
&=\sum_{k=1}^K \overline{\alpha}_{n_k}^k(t)^2
-\left(\sum_{k=1}^K \sqrt{\frac{n_k}{\boldsymbol{|n|}}}\,
\overline{\alpha}_{n_k}^k(t)\right)^{\!\!2}.
\end{align*}
Next, setting $\mathbb{B}_{K,\boldsymbol{n}}^*$ as in Theorem~\ref{theoremK},
we have
\begin{equation}\label{ecartK}
\overline{\boldsymbol{\xi}}_{K,\boldsymbol{n}}(t)-\mathbb{B}_{K,\boldsymbol{n}}^*(t)
=\Delta_{1,\boldsymbol{n}}(t)-\Delta_{2,\boldsymbol{n}}(t)
\end{equation}
where we put, for any $t\in\mathbb{R}$ and any
$\boldsymbol{n}=(n_1,\dots,n_K)\in\mathbb{N}^*$,
$$
\Delta_{1,\boldsymbol{n}}(t)=\sum_{k=1}^K \left(\overline{\alpha}_{n_k}^k(t)^2
-F_0(t)^2\,\mathbb{B}_{n_k}^k\!\big(F_0(t)\big)^2\right)\!,
$$
$$
\Delta_{2,\boldsymbol{n}}(t)=\left(\sum_{k=1}^K \sqrt{\frac{n_k}{\boldsymbol{|n|}}}\,
\overline{\alpha}_{n_k}^k(t)\right)^{\!\!2}
-\left(F_0(t) \sum_{k=1}^K\sqrt{\frac{n_k}{\boldsymbol{|n|}}}\,
\mathbb{B}_{n_k}^k\!\big(F_0(t)\big)\right)^{\!\!2}\!.
$$
By setting, for any $k\in\{1,\dots,K\}$, any $t\in\mathbb{R}$ and any
$\boldsymbol{n}=(n_1,\dots,n_K)\in\mathbb{N}^*$,
$$
\delta_{k,\boldsymbol{n}}(t)=\overline{\alpha}_{n_k}^k(t)
-F_0(t)\,\mathbb{B}_{n_k}^k\!\big(F_0(t)\big)
\quad\text{and}\quad\epsilon_{k,\boldsymbol{n}}(t)=
\overline{\alpha}_{n_k}^k(t)+F_0(t)\,\mathbb{B}_{n_k}^k\!\big(F_0(t)\big)
$$
and writing $\Delta_{1,\boldsymbol{n}}(t)$ and $\Delta_{2,\boldsymbol{n}}(t)$ as
$$
\Delta_{1,\boldsymbol{n}}(t)=\sum_{k=1}^K \delta_{k,\boldsymbol{n}}(t)
\,\epsilon_{k,\boldsymbol{n}}(t)\quad\text{and}\quad
\Delta_{2,\boldsymbol{n}}(t)=\sum_{k=1}^K \sqrt{\frac{n_k}{\boldsymbol{|n|}}}
\,\delta_{k,\boldsymbol{n}}(t) \,\sum_{k=1}^K\sqrt{\frac{n_k}{\boldsymbol{|n|}}}
\,\epsilon_{k,\boldsymbol{n}}(t),
$$
we derive the following inequalities:
\begin{align}
\sup_{t\in\mathbb{R}} |\Delta_{1,\boldsymbol{n}}(t)|
&
\leq\sum_{k=1}^K \left(\sup_{t\in\mathbb{R}} \left|\delta_{k,\boldsymbol{n}}(t)\right|\right)\!
\left(\sup_{t\in\mathbb{R}} \left|\epsilon_{k,\boldsymbol{n}}(t)\right|\right)\!,
\nonumber\\[-1ex]
\label{delta1}\\[-1ex]
\sup_{t\in\mathbb{R}} \big|\Delta_{2,\boldsymbol{n}}(t)\big|
&
\leq\left(\sum_{k=1}^K \sup_{t\in\mathbb{R}} \left|\delta_{k,\boldsymbol{n}}(t)\right|\right)\!
\left(\sum_{k=1}^K \sup_{t\in\mathbb{R}} \left|\epsilon_{k,\boldsymbol{n}}(t)\right|\right)\!.
\nonumber
\end{align}
By (\ref{corollaryapprox}), (\ref{lli}) and (\ref{estimbridge}),
we get the bounds, a.s., for each $k\in\{1,\dots,K\}$, as $n_k\to\infty$,
\begin{equation}\label{deltaeps}
\sup_{t\in\mathbb{R}} \left|\delta_{k,\boldsymbol{n}}(t)\right|
=O\!\left(\frac{\log n_k}{\sqrt{n_k}}\right)
\quad\text{and}\quad\sup_{t\in\mathbb{R}} \left|\epsilon_{k,\boldsymbol{n}}(t)\right|
=O\!\left(\!\sqrt{\log\log n_k}\,\right)\!.
\end{equation}
Finally, by putting (\ref{deltaeps}) into (\ref{delta1}),
and next into (\ref{ecartK}), we finish the proof of Theorem~\ref{theoremK}.
\hfill$\Box$

\subsection{Proof of Theorem~\ref{theorem1}}

In the computations below, the superscript ``$-$''
in the quantities $\mathbb{F}$, $\overline{\mathbb{F}}$, $\mathbb{U}$ and $\alpha$ refers
to the first $k$ observations while the superscript ``$+$'' refers to the last $n-k$ observations.
By (\ref{beta21}), we have the following representation for $\widetilde{\alpha}_n(s,t)$:
with probability $1$, for $n\in \mathbb{N}^*$, $t\in \mathbb{R}$ and $s\in(0,1)$,
\begin{align*}
\widetilde{\alpha}_n(s,t)=&\,\frac{\lfloor ns\rfloor(n-\lfloor ns\rfloor)}{n^{3/2}}
\left[ \left(\overline{\mathbb{F}}_{\lfloor ns\rfloor}^-(t)-\overline{F}(t)\right)
-\left( \overline{\mathbb{F}}_{n-\lfloor ns\rfloor}^+(t)-\overline{F}(t)\right) \right]
\\
=&\,\frac{\lfloor ns\rfloor(n-\lfloor ns\rfloor)}{2n^{3/2}}\left[
\left(\mathbb{F}_{\lfloor ns\rfloor}^-(t)^2-F(t)^2\right)
-\left( \mathbb{F}_{n-\lfloor ns\rfloor}^+(t)^2-F(t)^2\right) \right]
\\
&+\frac{\lfloor ns\rfloor(n-\lfloor ns\rfloor)}{2n^{3/2}}
\left( \frac{1}{\lfloor ns\rfloor}\,\mathbb{F}_{\lfloor ns\rfloor}^-(t)
-\frac{1}{n-\lfloor ns\rfloor}\,\mathbb{F}_{n-\lfloor ns\rfloor}^+(t)\right)\!.
\end{align*}
By using again the elementary identity $a^2 - b^2 = 2b(a-b) + (a -b)^2$, we can write
$\widetilde{\alpha}_n(s,t)$, a.s., for $n\in \mathbb{N}^*$,
$t\in\mathbb{ R}$ and $s\in (0, 1)$, in the form
\begin{equation}\label{equation78}
\widetilde{\alpha}_n(s,t)={\rm I}_n(s,t)-{\rm II}_n(s,t)+{\rm III}_n(s,t)+{\rm IV}_n(s,t)
\end{equation}
where
\begin{align*}
{\rm I}_n(s,t)&=\frac{\lfloor ns\rfloor(n-\lfloor ns\rfloor)}{n^{3/2}}\,F(t)
\left(\mathbb{F}_{\lfloor ns\rfloor}^-(t)-F(t)\right)
=\frac{\sqrt{\lfloor ns\rfloor}\,(n-\lfloor ns\rfloor)}{n^{3/2}}
\,F(t)\,\alpha_{\lfloor ns\rfloor}^-(t),
\\
{\rm II}_n(s,t)&=\frac{\lfloor ns\rfloor(n-\lfloor ns\rfloor)}{n^{3/2}}\,F(t)
\left(\mathbb{F}_{n-\lfloor ns\rfloor}^+(t)-F(t)\right)
=\frac{\lfloor ns\rfloor\sqrt{n-\lfloor ns\rfloor}}{n^{3/2}}
\,F(t)\,\alpha_{n-\lfloor ns\rfloor}^+(t),
\\
{\rm III}_n(s,t)&=\frac{\lfloor ns\rfloor(n-\lfloor ns\rfloor)}{2n^{3/2}}
\left[\left(\mathbb{F}_{\lfloor ns\rfloor}^-(t)-F(t)\right)^{\!2}
-\left(\mathbb{F}_{n-\lfloor ns\rfloor}^+(t)-F(t)\right)^{\!2}\right]
\\
&=\frac{1}{2n^{3/2}}\left((n-\lfloor ns\rfloor)\,\alpha_{\lfloor ns\rfloor}^-(t)^2
-\lfloor ns\rfloor\alpha_{n-\lfloor ns\rfloor}^+(t)^2\right)\!,
\\
{\rm IV}_n(s,t)&=\frac{1}{2n^{3/2}}\left((n-\lfloor ns\rfloor)\,
\mathbb{F}_{\lfloor ns\rfloor}^-(t)-
\lfloor ns\rfloor\mathbb{F}_{n-\lfloor ns\rfloor}^+(t)\right)\!.
\end{align*}
By (\ref{empililaeae}) and by using the fact that $\mathbb{F}_n$ is bounded, we get that,
with probability $1$, as $n \to\infty$, uniformly in $s$ and~$t$,
$$
\mathrm{III}_n(s,t)=O\!\left(\frac{\log \log n}{\sqrt{n}}\right)
\quad\text{and}\quad
\mathrm{IV}_n(s,t)=O\!\left(\frac{1}{\sqrt{n}}\right)\!.
$$
which, inserted in (\ref{equation78}), provides,
with probability $1$, as $n \to\infty$, uniformly in $s$ and $t$,
\begin{equation}\label{equation78bis}
\widetilde{\alpha}_n(s,t)={\rm I}_n(s,t)-{\rm II}_n(s,t)
+O\!\left(\frac{\log \log n}{\sqrt{n}}\right)\!.
\end{equation}
Let us introduce, for $n\in \mathbb{N}^*$ and $s,u\in(0,1),$
\begin{align*}
{\rm I}_n'(s,u)&=\frac{\lfloor ns\rfloor(n-\lfloor ns\rfloor)}{n^{3/2}}
\,u\left(\mathbb{U}_{\lfloor ns\rfloor}^-(u)-u\right)
=\frac{n-\lfloor ns\rfloor}{n^{3/2}} \,u\,\Bigg(\sum_{i=1}^{\lfloor ns\rfloor}
\left(\mathbbm{1}_{\{ U_i\leq u\}}-u\right)\!\Bigg)\!,
\\
{\rm II}_n'(s,u)&=\frac{\lfloor ns\rfloor(n-\lfloor ns\rfloor)}{n^{3/2}}
\,u\left(\mathbb{U}_{n-\lfloor ns\rfloor}^+(u)-u\right)
=\frac{\lfloor ns\rfloor}{n^{3/2}} \,u \,\Bigg(\sum_{i=\lfloor ns\rfloor+1}^n
\left(\mathbbm{1}_{\{ U_i\leq u\}}-u\right)\!\Bigg)\!,
\end{align*}
and
$$
\delta_n(s,u)={\rm I}_n'(s,u)-{\rm II}_n'(s,u).
$$
Then, we have the following equalities
$$
{\rm I}_n'(s,F(t))={\rm I}_n(s,t),\quad {\rm II}_n'(s,F(t))={\rm II}_n(s,t),
$$
and (\ref{equation78bis}) becomes, with probability $1$, as $n \to\infty$,
uniformly in $s$ and $t$,
\begin{equation}\label{equatioreffree}
\widetilde{\alpha}_n(s,t)= \delta_n(s,F(t))+O\!\left(\frac{\log \log n}{\sqrt{n}}\right)\!.
\end{equation}
Now, observe that
\begin{align}\label{h-n15,31}
\delta_n(s,u)&=\frac{u}{\sqrt{n}}\Bigg( \sum_{i=1}^{\lfloor ns\rfloor}
\left(\mathbbm{1}_{\{ U_i\leq u\}}-u\right)-\frac{\lfloor ns\rfloor}{n}
\sum_{i=1}^n\left(\mathbbm{1}_{\{ U_i\leq u\}}-u\right)\Bigg)
\\[1ex]
\label{h-n15,31bis}
&=\frac{u}{\sqrt{n}}\Bigg( \frac{n-\lfloor ns\rfloor}{n}
\sum_{i=1}^n\left(\mathbbm{1}_{\{ U_i\leq u\}}-u\right)-\sum_{i=\lfloor ns\rfloor+1}^n
\left(\mathbbm{1}_{\{ U_i\leq u\}}-u\right)\Bigg).
\end{align}
We know from \cite{KMT1975} and \cite{Csorgo1997} that, almost surely, as
$n\to\infty$,
\begin{align}\label{h-n15,36}
\sup_{s\in[0,1/2]}\sup_{u\in[0,1]}\Bigg| \sum_{i=1}^{\lfloor ns\rfloor}
\left(\mathbbm{1}_{\{ U_i\leq u\}} -u\right)
-\mathbb{K}_2(\lfloor ns\rfloor,u)\Bigg| &=O\!\left((\log n)^2\right)\!,
\\[1ex]
\label{h-n15,36bis}
\sup_{s\in[1/2,1]}\sup_{u\in[0,1]}\Bigg| \sum_{i=\lfloor ns\rfloor+1}^n
\left(\mathbbm{1}_{\{ U_i\leq u\}} -u\right)
-\mathbb{K}_1(\lfloor ns\rfloor,u)\Bigg|&=O\!\left((\log n)^2\right)\!.
\end{align}
Notice that we have the following decomposition
$$
\sum_{i=1}^n\left( \mathbbm{1}_{\{ U_i\leq u\}} -u\right)
=\sum_{i=1}^{\lfloor n/2\rfloor}\left( \mathbbm{1}_{\{ U_i\leq u\}} -u\right)
+\sum_{i=\lfloor n/2\rfloor+1}^n\left( \mathbbm{1}_{\{ U_i\leq u\}} -u\right)\!.
$$
Hence, by adding (\ref{h-n15,36}) and (\ref{h-n15,36bis}), we readily infer that,
almost surely, as $n\to\infty$,
\begin{equation}\label{prth7}
\sup_{u\in[0,1]}\left| \sum_{i=1}^n\left(\mathbbm{1}_{\{ U_i\leq u\}} -u\right)
-\big( \mathbb{K}_1(\lfloor n/2\rfloor,u)
+\mathbb{K}_2(\lfloor n/2\rfloor,u)\big)\right|=O\!\left((\log n)^2\right)\!.
\end{equation}
As a byproduct, from (\ref{h-n15,31})--(\ref{prth7})
and recalling the definition of $\overline{\mathbb{K}}_n$ given just before
Theorem~\ref{theorem1}, we deduce that, almost surely, as $n\to\infty$,
\begin{equation}\label{prth7bis}
\sup_{s\in[0,1]}\sup_{u\in[0,1]} \left|\delta_n(s,u)-\overline{\mathbb{K}}_n(s,u)\right|
=O\!\left(\frac{(\log n)^2}{\sqrt{n}}\right)\!.
\end{equation}
We finally conclude from (\ref{equatioreffree}) and (\ref{prth7bis}) by using
the triangle inequality: almost surely, as $n\to\infty$,
\begin{align*}
\sup_{s\in[0,1]}\sup_{t\in \mathbb{R}} \left|\widetilde{\alpha}_n(s,t)
-\overline{\mathbb{K}}_n(s,F(t))\right|
&
\leq \sup_{s\in[0,1]}\sup_{t\in \mathbb{R}} \big|\widetilde{\alpha}_n(s,t)
-\delta_n(s,F(t))\big|+\sup_{s\in[0,1]}\sup_{u\in[0,1]}
\left|\delta_n(s,u)-\overline{\mathbb{K}}_n(s,u)\right|
\\
&
= O\!\left(\frac{(\log n)^2}{\sqrt{n}}\right)
+O\!\left(\frac{\log \log n}{\sqrt{n}}\right)
=O\!\left(\frac{(\log n)^2}{\sqrt{n}}\right)\!.
\end{align*}
This completes the proof of Theorem~\ref{theorem1}.
\hfill$\Box$\\

\subsection{Proof of Corollary~\ref{coroltau}}

Straightforward algebra yields, for any $s,t,u,v\in[0,1]$,
$$
\mathbb{E}\!\left(\overset{\text{\tiny o}}{\mathbb{K}}_n(s,u)\,
\overset{\text{\tiny o}}{\mathbb{K}}_n(t,v)\right)=\frac{1}{n}\,(u\wedge v-uv)\,\psi_n(s,t)
$$
with
$$
\psi_n(s,t)=\begin{cases}
\lfloor n(s\wedge t)\rfloor-s\lfloor nt\rfloor-t\lfloor ns\rfloor+2\lfloor n/2\rfloor st
&\text{for $s,t\in[0,1/2]$,}
\\[1ex]
\lfloor n(1-s\vee t)\rfloor-(1-s)\lfloor n(1-t)\rfloor
&
\\[-0.5ex]
\quad-(1-t)\lfloor n(1-s)\rfloor+2\lfloor n/2\rfloor (1-s)(1-t)
&\text{for $s,t\in[1/2,1]$,}
\\[1ex]
s\lfloor n(1-t)\rfloor+(1-t)\lfloor ns\rfloor-2\lfloor n/2\rfloor s(1-t)
&\text{for $s\in[0,1/2],\,t\in[1/2,1]$,}
\\[1ex]
(1-s)\lfloor nt\rfloor+t\lfloor n(1-s)\rfloor-2\lfloor n/2\rfloor (1-s)t
&\text{for $s\in[1/2,1],\,t\in[0,1/2]$.}
\end{cases}
$$
We immediately see that
$$
\lim_{n\to\infty}\frac{1}{n}\,\psi_n(s,t)=s\wedge t-st
$$
and then
$$
\lim_{n\to\infty}\mathbb{E}\!\left(\overset{\text{\tiny o}}{\mathbb{K}}_n(t,v)
\overset{\text{\tiny o}}{\mathbb{K}}_n(s,u)\right)=(s\wedge t-st)\,(u\wedge v-uv)
=\mathbb{E}\!\left(\overset{\text{\tiny o}}{\mathbb{K}}(t,v)
\overset{\text{\tiny o}}{\mathbb{K}}(s,u)\right)
$$
where $\overset{\text{\tiny o}}{\mathbb{K}}$ is the tied-down Kiefer process
on $[0,1]\!\times\![0,1]$.
This proves the convergence of Gaussian processes in distribution, as $n\to\infty$,
$$
\overset{\text{\tiny o}}{\mathbb{K}}_n\stackrel{\mathcal{L}}{\longrightarrow}
\overset{\text{\tiny o}}{\mathbb{K}},
$$
which in turn, together with Theorem~\ref{theorem1}, entails Corollary
\ref{coroltau}.

\subsection{Proof of Theorem~\ref{theoremBurkeetal}}

Let us write $\overline{\widehat{\alpha}}_n(t)$ as follows: with probability~$1$,
for any $t\in\mathbb{R},\,n\in\mathbb{N}^*$,
\begin{align*}
\overline{\widehat{\alpha}}_n(t)&=\sqrt{n}\left( \overline{\mathbb{F}}_n(t)
-\overline{F}\big(t,\widehat{\boldsymbol{\theta}}_n\big)\right)
=\frac{\sqrt{n}}{2}\left(\mathbb{F}_n(t)^2-F\big(t,\widehat{\boldsymbol{\theta}}_n\big)^2\right)
+\frac{1}{2\sqrt{n}}\,\mathbb{F}_n(t)
\\
&=\sqrt{n}\,F\big(t,\widehat{\boldsymbol{\theta}}_n\big)\left(\mathbb{F}_n(t)
-F\big(t,\widehat{\boldsymbol{\theta}}_n\big)\right)+\frac{\sqrt{n}}{2}\left(\mathbb{F}_n(t)
-F\big(t,\widehat{\boldsymbol{\theta}}_n\big)\right)^2+\frac{1}{2\sqrt{n}}\,\mathbb{F}_n(t)
\\
&=F\big(t,\widehat{\boldsymbol{\theta}}_n\big)\,\widehat{\alpha}_n(t)
+\frac{1}{2\sqrt{n}}\,\widehat{\alpha}_n(t)^2+\frac{1}{2\sqrt{n}}\,\mathbb{F}_n(t).
\end{align*}
Straightforward algebra yields, a.s., for any $t\in\mathbb{R},\,n\in\mathbb{N}^*$,
\begin{align}
\overline{\widehat{\alpha}}_n(t)-\overline{G}_n(t)=&\,\left(F\big(t,\widehat{\boldsymbol{\theta}}_n\big)
+\frac{1}{\sqrt{n}}\,G_n(t)\right) \big(\widehat{\alpha}_n(t)-G_n(t)\big)
+\frac{1}{2\sqrt{n}}\,\big(\widehat{\alpha}_n(t)-G_n(t)\big)^2
\nonumber\\
&\,+\left(F\big(t,\widehat{\boldsymbol{\theta}}_n\big)-F(t,\boldsymbol{\theta}_0)\right)G_n(t)
+\frac{1}{2\sqrt{n}}\,G_n(t)^2+\frac{1}{2\sqrt{n}}\,\mathbb{F}_n(t).
\label{ecartalphahat}
\end{align}
Recall the definition of $\overline{\boldsymbol{\varepsilon}}_n$ given
just before Theorem~\ref{theoremBurkeetal} and set
$$
\boldsymbol{\eta}_n:=\sup_{t\in \mathbb{R}}\left|\widehat{\alpha}_n(t)-G_n(t)\right|\!.
$$
Using again the inequalities $0\leq F(t)\leq 1$ and $0\leq \mathbb{F}_n(t)\leq 1$
and roughly bounding $1/2$ by $1$, we extract from (\ref{ecartalphahat}) that,
a.s., for any $t\in\mathbb{R},\,n\in\mathbb{N}^*$,
\begin{align}
\overline{\boldsymbol{\varepsilon}}_n
\leq&\,
\left(1+\frac{1}{\sqrt{n}}
\sup_{t\in \mathbb{R}}|G_n(t)|\right)\boldsymbol{\eta}_n+\frac{1}{\sqrt{n}}\,
\boldsymbol{\eta}_n^2 +\frac{1}{\sqrt{n}}\left(1+\sup_{t\in \mathbb{R}}G_n(t)^2\right)
\nonumber\\
&
\;+\sup_{t\in \mathbb{R}}\left|F\big(t,\widehat{\boldsymbol{\theta}}_n\big)
-F(t,\boldsymbol{\theta}_0)\right| \,\sup_{t\in \mathbb{R}}|G_n(t)|.
\label{reffequa}
\end{align}
We know from Theorem 3.1 of \cite{Burke-Csorgo} that $\boldsymbol{\eta}_n$
satisfies the same limiting results that those displayed in
Theorem~\ref{theoremBurkeetal} for $\overline{\boldsymbol{\varepsilon}}_n$.
Next, we need to derive some bounds for $\sup_{t\in \mathbb{R}}|G_n(t)|$ and
$\sup_{t\in \mathbb{R}}\big|F\big(t,\widehat{\boldsymbol{\theta}}_n\big)
-F(t,\boldsymbol{\theta}_0)\big|$ as $n\to\infty$. First, by
(\ref{equation11KIE}), we have, almost surely, as $n\to\infty$,
$$
\sup_{t\in \mathbb{R}}\left|\alpha_n(t)-\frac{1}{\sqrt{n}}\,\mathbb{K}(n,F(t))\right|
=O\!\left(\frac{(\log n)^2}{\sqrt{n}}\right)
$$
from which we deduce, due to (\ref{empililaeae}), almost surely, as $n\to\infty$,
$$
\sup_{u\in[0,1]}|\mathbb{K}(n,u)|=O\!\left(\!\sqrt{n \log \log n}\,\right)\!.
$$
Notice that the same bound holds for the Brownian motion $\{\mathbf{W}(n):
n\in\mathbb{N}^*\}$ introduce before Theorem~\ref{theoremBurkeetal}. Hence,
by Condition~(iv) and the definition of $G_n(t)$, with probability $1$,
as $n\to\infty$,
\begin{equation}\label{reffffreco}
\sup_{t\in \mathbb{R}}|G_n(t)|=O\!\left(\!\sqrt{n \log \log n}\,\right)\!.
\end{equation}
On the other hand, using the one-term Taylor expansion of $F(\cdot,\boldsymbol{\theta})$ with respect
to $\boldsymbol{\theta}_0$, there exists $\boldsymbol{\theta}_n^*$ lying in the segment
$\big[\boldsymbol{\theta}_0, \widehat{\boldsymbol{\theta}}_n\big]$ such that
\begin{equation}\label{reffffreco23}
F\big(t,\widehat{\boldsymbol{\theta}}_n\big)-F(t,\boldsymbol{\theta}_0)
=\left(\widehat{\boldsymbol{\theta}}_n-\boldsymbol{\theta}_0\right)
\nabla_{\boldsymbol{\theta}}F(t,\boldsymbol{\theta}_n^*)^{\top}.
\end{equation}

In case (a) of Theorem~\ref{theoremBurkeetal},
$\sqrt{n}\,\big(\widehat{\boldsymbol{\theta}}_n-\boldsymbol{\theta}_0\big)$
is asymptotically normal and then
$n^{1/4}\big(\widehat{\boldsymbol{\theta}}_n-\boldsymbol{\theta}_0\big)$
tends to zero as $n \to\infty$ in probability.
Therefore, by (\ref{reffffreco}) and (\ref{reffffreco23}),
$\sup_{t\in \mathbb{R}}\big|F\big(t,\widehat{\boldsymbol{\theta}}_n\big)
-F(t,\boldsymbol{\theta}_0)\big|\sup_{t\in \mathbb{R}}|G_n(t)|$
also tends to zero as $n\to\infty$ in probability. Putting this into (\ref{reffequa}),
we easily complete the proof of Theorem~\ref{theoremBurkeetal} in this case. In cases (b) and (c) of Theorem~\ref{theoremBurkeetal}, referring to \cite{Burke-Csorgo} p.~779, we have the
following bound for $\widehat{\boldsymbol{\theta}}_n-\boldsymbol{\theta}_0$:
almost surely, as $n\to\infty$,
$$
\widehat{\boldsymbol{\theta}}_n-\boldsymbol{\theta}_0
=O\!\left(\!\sqrt{\frac{\log \log n}{n}}\,\right)\!.
$$
By putting this into (\ref{reffffreco23}) and next in (\ref{reffequa}) with the aid of
(\ref{reffffreco}), we complete the proof of Theorem~\ref{theoremBurkeetal} in these two cases.

Finally, concerning $\overline{\widehat{G}}_n(t)$, we have
$$
\overline{\widehat{G}}_n(t)-\overline{G}_n(t)
=F\big(t,\widehat{\boldsymbol{\theta}}_n\big)\left(\widehat{G}_n(t)-G_n(t)\right)
+\left(F\big(t,\widehat{\boldsymbol{\theta}}_n\big)-F(t,{\boldsymbol{\theta}}_0)\right)G_n(t),
$$
from which we deduce
$$
\sup_{t\in \mathbb{R}}\left|\overline{\widehat{G}}_n(t)-\overline{G}_n(t)\right|\leq
\sup_{t\in \mathbb{R}}\left|\widehat{G}_n(t)-G_n(t)\right|
+\sup_{t\in \mathbb{R}}\left|F\big(t,\widehat{\boldsymbol{\theta}}_n\big)
-F(t,{\boldsymbol{\theta}}_0)\right|\sup_{t\in \mathbb{R}}|G_n(t)|.
$$
Using the same bounds than previously, we immediately derive (\ref{eqreferrrz11}).\hfill$\Box$\\

\subsection{Proof of Theorem~\ref{Res1}}

In a similar way as in \cite{BassKhoshnevisan1993} (refer also to \cite{BassKhoshnevisan1992}),
we consider the following normalized local time
\begin{equation}\label{BaKh0}
\Lambda_n(x,t):=\frac{1}{\sqrt{n}}\,\lambda\!\left(\sqrt{n}\,x,\lfloor nt\rfloor\right)
\quad\text{for}\quad x\in \mathbb{R},\, t\in[0,1],\,n\in \mathbb{N}^*,
\end{equation}
where the local time $\lambda$ is defined by (\ref{localtime}).
We first decompose $\Lambda_n(x,t)^2$ into the sum of two components, by writing
\begin{align*}
\Lambda_n(x,t)^2&=\frac{1}{n}\sum_{i=1}^{\lfloor nt\rfloor}
\sum_{j=1}^{\lfloor nt\rfloor} \mathbbm{1}_I\!\left(S_i-\sqrt{n}\,x\right)\mathbbm{1}_I\!\left(S_j-\sqrt{n}\,x\right)
\\
&=\frac{1}{n}\sum_{i=1}^{\lfloor nt\rfloor} \mathbbm{1}_I\!\left(S_i-\sqrt{n}\,x\right)
+\frac{2}{n}\sum_{1\leq i<j \leq \lfloor nt\rfloor }
\mathbbm{1}_I\!\left(S_i-\sqrt{n}\,x\right)\mathbbm{1}_I\!\left(S_j-\sqrt{n}\,x\right)\!.
\end{align*}
This, in turn, implies that
$$
\int_{\mathbb{R}}\Lambda_n(x,t)^2\,dx=\frac{1}{n}\int_{\mathbb{R}}
\sum_{i=1}^{\lfloor nt\rfloor} \mathbbm{1}_I\!\left(S_i-\sqrt{n}\,x\right)dx
+\frac{2}{n}\sum_{1\leq i<j \leq \lfloor nt\rfloor }\int_{\mathbb{R}}
\mathbbm{1}_I\!\left(S_i-\sqrt{n}\,x\right)\mathbbm{1}_I\!\left(S_j-\sqrt{n}\,x\right)dx.
$$
from which, by observing that
$\int_{\mathbb{R}}\mathbbm{1}_I\!\left(S_i-\sqrt{n}\,x\right)dx=1/\sqrt{n}$
and that
$$
\sum_{1\leq i<j \leq \lfloor nt\rfloor }\int_{\mathbb{R}}
\mathbbm{1}_I\!\left(S_i-\sqrt{n}\,x\right)\mathbbm{1}_I\!\left(S_j-\sqrt{n}\,x\right)dx
=\frac{1}{\sqrt{n}} \,L_n(t),
$$
we readily infer that
\begin{equation}\label{HuKh0}
L_n(t)=\frac{1}{2}\,n^{3/2}\int_{\mathbb{R }}\Lambda_n(x,t)^2\,dx
-\frac{1}{2}\lfloor nt\rfloor.
\end{equation}
In order to evaluate the right hand-side of (\ref{HuKh0}), first, remark that
\begin{equation}\label{ecartl}
\left|\int_{\mathbb{R}}\Lambda_n(x,t)^2\,dx
-\int_{\mathbb{R}}l_n(x,t)^2\,dx\right|\leq \sup_{x\in \mathbb{R}}|\Lambda_n(x,t)-l_n(x,t)|
\int_{\mathbb{R}}(\Lambda_n(x,t)+l_n(x,t))\,dx.
\end{equation}
By combining (\ref{BaKh0}) with (\ref{BaKh1}), we obtain, with probability $1$, as $n\to\infty$,
\begin{equation}\label{ecartlbis}
\sup_{x\in \mathbb{R}}\sup_{t\in[0,1]}|\Lambda_n(x,t)-l_n(x,t)|
=O\!\left(n^{-1/4}(\log n)^{1/2}(\log \log n)^{1/4}\right)\!.
\end{equation}
Moreover, using the following relations satisfied by the normalized local times:
$$
\int_{\mathbb{R}}\Lambda_n(x,t)\,dx=\int_{\mathbb{R}}l_n(x,t)\,dx=\lfloor nt\rfloor/n,
$$
we get, for any $t\in[0,1]$ and any $n\in \mathbb{N}^*$,
\begin{equation}\label{ecartlter}
\int_{\mathbb{R}}(\Lambda_n(x,t)+l_n(x,t))\,dx\leq 2.
\end{equation}
Then, putting (\ref{ecartlbis}) and (\ref{ecartlter}) into (\ref{ecartl})
yields, with probability $1$, as $n\to\infty$,
\begin{equation}\label{ecartlretour}
\sup_{t\in[0,1]}\left|\int_{\mathbb{R}}\Lambda_n(x,t)^2\,dx
-\int_{\mathbb{R}}l_n(x,t)^2\,dx\right|
=O\!\left(n^{-1/4}(\log n)^{1/2}(\log \log n)^{1/4}\right)\!.
\end{equation}
Finally, from (\ref{HuKh0}) and (\ref{ecartlretour}), we infer, with
probability~$1$, as $n\to\infty$, uniformly in $t$,
$$
L_n(t)=\frac{1}{2}\,n^{3/2}\int_{\mathbb{R}} l_n(x,t)^2\,dx
+O\!\left(n^{5/4}(\log n)^{1/2}(\log \log n)^{1/4}\right)\!.
$$
This completes the proof of Theorem~\ref{Res1}.
\hfill$\Box$

\subsection{Proof of Corollary~\ref{corLn}}

We imitate the proof of Lemma~2.2 of \cite{HuKhoshnevisan2010}
yielding $\int_{\mathbb{R}}l(x,n)^2\,dx=n^{3/2+o(1)}$.
For this, let us observe that
$$
\int_{\mathbb{R}} l(x,n)^2\,dx=\int_{I_n} l(x,n)^2\,dx
$$
where $I_n$ is the range of the Brownian motion over the time interval $[0,n]$, i.e.,
$$
I_n=\Big[\inf_{s\in[0,n]}\mathbb{W}(s),\sup_{s\in[0,n]}\mathbb{W}(s)\Big]
$$
and let us introduce the corresponding oscillation
$$
\mathbb{Z}_n:=\sup_{s\in[0,n]}\mathbb{W}(s)-\inf_{s\in[0,n]}\mathbb{W}(s).
$$
By using Cauchy-Schwarz inequality, we write the inequalities
$$
\frac{1}{\mathbb{Z}_n}\left(\int_{\mathbb{R}} l(x,n)\,dx\right)^{\!2}
\leq \int_{\mathbb{R}} l(x,n)^2\,dx
\leq \int_{\mathbb{R}} l(x,n)\,dx\times \sup_{x\in\mathbb{R}} l(x,n)
$$
which can be rewritten, since $\int_{\mathbb{R}}l(x,n)\,dx=n$, as
\begin{equation}\label{inequalitylnbis}
\frac{n^2}{\mathbb{Z}_n} \leq \int_{\mathbb{R}} l(x,n)^2\,dx
\leq n\sup_{x\in\mathbb{R}} l(x,n).
\end{equation}
Now, by invoking Paul L\'evy's (\cite{Csorgo1981REVESZ} p.~36) and
\cite{Kesten1965}'s laws of the iterated logarithm which respectively state that
$$
\limsup_{\tau\to\infty} \frac{\sup_{s\in[0,\tau]}|\mathbb{W}(s)|}{\sqrt{\tau\log\log \tau}}
=\sqrt{2}\quad\text{and}\quad
\limsup_{\tau\to\infty} \frac{l(x,\tau)}{\sqrt{\tau\log\log \tau}}=\sqrt{2}\quad\text{a.s.},
$$
we obtain that, with probability $1$, there exist
two positive constants $\kappa_1'$ and $\kappa_2'$ such that, for large enough $n$,
\begin{equation}\label{lliZ}
\mathbb{Z}_n\leq \frac{1}{\kappa_1'}\sqrt{n\log\log n}
\quad\text{and}\quad \sup_{x\in\mathbb{R}} l(x,n)\leq\kappa_2'\sqrt{n\log\log n}.
\end{equation}
By putting (\ref{lliZ}) into (\ref{inequalitylnbis}), we deduce (\ref{inequalityln}).
Finally, by observing that
$\int_{\mathbb{R}} l_n(x,t)^2\,dx={n^{-3/2}}\int_{\mathbb{R}} l(x,\lfloor nt\rfloor)^2\,dx$
and appealing to Theorem~\ref{Res1}, we get (\ref{inequalityLn}).
The proof of Corollary~\ref{corLn} is completed.

\section*{Acknowledgement}
The authors are grateful to the Editor, an Associate editor, and the referees for thorough proofreading,
numerous comments and several inspiring questions which led to a considerable improvement of the
presentation.


\end{document}